\documentclass[titlepage,draft,12pt]{article} 
\usepackage{amsfonts,latexsym} 
\usepackage[a4paper]{geometry}

\date{}

\newcommand{\ep}{\varepsilon}
\newcommand{\qed}{{\penalty 10000\mbox{$\quad\Box$}}}
\newcommand{\re}{\mathbb{R}}

\newcommand{\n}{\mathbb{N}}
\newcommand{\cep}{c_{\ep}}
\newcommand{\rep}{r_{\ep}}
\newcommand{\roep}{\rho_{\ep}}
\newcommand{\uepl}{u_{\ep,\ell}}
\newcommand{\roepl}{\rho_{\ep,\ell}}
\newcommand{\ueph}{u_{\ep,h}}
\newcommand{\roeph}{\rho_{\ep,h}}
\newcommand{\uephh}{u_{\ep,h,\mu}}
\newcommand{\roephh}{\rho_{\ep,h,\mu}}
\newcommand{\ueplm}{u_{\ep,\ell,\mu}}
\newcommand{\roeplm}{\rho_{\ep,\ell,\mu}}
\newcommand{\teteplm}{\theta_{\ep,\ell,\mu}}
\newcommand{\uep}{u_{\ep}}
\newcommand{\ul}{u_{\ell}}
\newcommand{\uh}{u_{h}}
\newcommand{\gep}{g_{\ep}}
\newcommand{\Dep}{\mathcal{D}_{\ep}}
\newcommand{\Eep}{\mathcal{E}_{\ep}}
\newcommand{\Gep}{\mathcal{G}_{\ep}}
\newcommand{\tetep}{\theta_{\ep}}
\newcommand{\auq}[1]{|A^{1/2}#1|^{2}}
\newcommand{\auqg}[1]{|A^{1/2}#1|^{2\gamma}}
\newcommand{\dat}{D(A^{3/2})}
\newcommand{\da}{D(A)}
\newcommand{\dau}{D(A^{1/2})}

\newcommand{\dg}{\eta}

\newtheorem{thm}{Theorem}[section]
\newtheorem{thmbibl}{Theorem}
\newtheorem{propbibl}[thmbibl]{Proposition}
\newtheorem{rmk}[thm]{Remark}
\newtheorem{prop}[thm]{Proposition}
\newtheorem{defn}[thm]{Definition}

\newtheorem{ex}[thm]{Example}
\newtheorem{lemma}[thm]{Lemma}

\title{Hyperbolic-parabolic singular perturbation for mildly
degenerate Kirchhoff equations: decay-error estimates}

\author{Marina Ghisi\vspace{1ex}\\ 
{\normalsize Universit\`a degli Studi di Pisa} \\
{\normalsize Dipartimento di Matematica ``Leonida Tonelli''}\\ 
{\normalsize PISA (Italy)}\\
{\normalsize e-mail: \texttt{ghisi@dm.unipi.it}}
\and
Massimo Gobbino\vspace{1ex}\\ 
{\normalsize Universit\`a degli Studi di Pisa} \\
{\normalsize Dipartimento di Matematica Applicata ``Ulisse Dini''}\\ 
{\normalsize PISA (Italy)}\\  
{\normalsize e-mail: \texttt{m.gobbino@dma.unipi.it}}}

\begin{document}
\maketitle
\begin{abstract}
	We consider \emph{degenerate} Kirchhoff equations with a small
	parameter $\ep$ in front of the second-order time-derivative.  It
	is well known that these equations admit global solutions when
	$\ep$ is small enough, and that these solutions decay as $t\to
	+\infty$ with the same rate of solutions of the limit problem
	(of parabolic type).
	
	In this paper we prove decay-error estimates for the difference
	between a solution of the hyperbolic problem and the solution of
	the corresponding parabolic problem.  These estimates show in the
	same time that the difference tends to zero both as $\ep\to
	0^{+}$, and as $t\to +\infty$.  Concerning the decay rates, it
	turns out that the difference decays \emph{faster} than the two
	terms separately (as $t\to +\infty$).
	
	Proofs involve a nonlinear step where we separate Fourier
	components with respect to the lowest frequency, followed by a
	linear step where we exploit weighted versions of classical
	energies.
	
\vspace{1cm}

\noindent{\bf Mathematics Subject Classification 2000 (MSC2000):}
35B25, 35L70, 35L80.

\vspace{1cm} 

\end{abstract}

 
\section{Introduction}

Let $H$ be a separable real Hilbert space.  For every $x$ and $y$ in
$H$, $|x|$ denotes the norm of $x$, and $\langle x,y\rangle$ denotes
the scalar product of $x$ and $y$.  Let $A$ be a self-adjoint linear
operator on $H$ with dense domain $D(A)$.  We assume that $A$ is
nonnegative, namely $\langle Ax,x\rangle\geq 0$ for every $x\in D(A)$,
so that for every $\alpha\geq 0$ the power $A^{\alpha}x$ is defined
provided that $x$ lies in a suitable domain $D(A^{\alpha})$.

We consider the Cauchy problem
\begin{equation}
	\ep\uep''(t)+\uep'(t)+\auqg{\uep(t)}A\uep(t)=0
	\quad\quad
	\forall t\geq 0,
	\label{pbm:h-eq}
\end{equation}
\begin{equation}
	\uep(0)=u_0,\hspace{3em}\uep'(0)=u_1,
	\label{pbm:h-data}
\end{equation}
where $\ep>0$ and $\gamma\geq 1$ are real parameters, and 
$(u_{0},u_{1})\in\da\times\dau$ are initial conditions satisfying the 
mild nondegeneracy condition
\begin{equation}
	A^{1/2}u_{0}\neq 0.
	\label{hp:mdg}
\end{equation}

The singular perturbation problem in its generality consists in
proving the convergence of solutions of (\ref{pbm:h-eq}),
(\ref{pbm:h-data}) to solutions of the first order problem
\begin{equation}
	u'(t)+ \auqg{u(t)}Au(t)=0
	\quad\quad
	\forall t\geq 0,
	\label{pbm:p-eq}
\end{equation}
\begin{equation}
	u(0)=u_{0},
	\label{pbm:p-data}
\end{equation}
obtained setting formally $\ep=0$ in (\ref{pbm:h-eq}), and
omitting the second initial condition in~(\ref{pbm:h-data}).
Following the approach introduced by \textsc{J.\ L.\
Lions}~\cite{lions} in the linear case, one defines the corrector
$\tetep(t)$ as the solution of the second order \emph{linear} problem
\begin{equation}
	\ep\tetep''(t)+\tetep'(t)=0 \hspace{2em}
	\forall t\geq 0,
	\label{pbm:tetep-eq}
\end{equation}
\begin{equation}
	\tetep(0)=0,\hspace{2em}\tetep'(0)=u_1+
	\auqg{u_{0}}Au_{0}=:w_{0}.
	\label{pbm:tetep-data}
\end{equation}

It is easy to see that $\tetep'(0)=\uep'(0)-u'(0)$, hence this
corrector keeps into account the boundary layer due to the loss of one
initial condition.  Finally, one defines $\rep(t)$ and $\roep(t)$ in
such a way that
$$\uep(t)=u(t)+\tetep(t)+\rep(t)=u(t)+\roep(t)\quad\quad\forall t\geq
0.$$

With these notations, the singular perturbation problem consists in
proving that $\rep(t)\to 0$ or $\roep(t)\to 0$ in some sense as
$\ep\to 0^{+}$. 

The singular perturbation problem for Kirchhoff equations has
generated a considerable literature in the last 30 years.  The state
of the art has been recently presented in the
survey~\cite{gg:survey-diss}, where more general nonlinearities and
more general dissipative terms have also been considered.
In~\cite{gg:survey-diss} the general problem has been split into six
subproblems, which we list below.
\begin{enumerate}
	\renewcommand{\labelenumi}{(P\arabic{enumi})} 
	\item Global existence and decay estimates for the parabolic
	problem.

	\item Local existence for the hyperbolic problem and
	local-in-time error estimates on $\roep(t)$ and $\rep(t)$.

	\item Global existence for the hyperbolic problem.
	
	\item Decay estimates for solutions of the hyperbolic problem (as
	$t\to +\infty$).

	\item Global-in-time error estimates for the singular perturbation
	problem, which means \emph{time-independent} estimates on
	$\roep(t)$ or $\rep(t)$ as $\ep\to 0^{+}$.

	\item Decay-error estimates for the singular perturbation 
	problem, which means estimates such as
	\begin{equation}
		|A^{\alpha}\roep(t)|\leq\omega(\ep)\sigma(t)
		\quad\mbox{or}\quad
		|A^{\alpha}\rep(t)|\leq\omega(\ep)\sigma(t),
		\label{est:de-model}
	\end{equation}
	where of course the convergence rate $\omega(\ep)$ tends to 0 as
	$\ep\to 0^{+}$, and the decay rate $\sigma(t)$ tends to 0 as
	$t\to+\infty$.  Decay-error estimates are the meeting point of
	subproblems (P4) and (P5), and they represent the ultimate goal of
	the theory.
\end{enumerate}

Subproblem~(P1) is well understood
(see~\cite{bw,bernstein,k-par,miletta}).  The result is that problem
(\ref{pbm:p-eq}), (\ref{pbm:p-data}) has a unique global solution for
every $u_{0}\in D(A)$ (and even for less regular data), and this
solution decays at infinity as solutions of the ordinary differential
equation
\begin{equation}
	y'+|y|^{2\gamma}y=0,
	\label{eqn:ODE-model}
\end{equation}
which is just the special case of (\ref{pbm:p-eq}) where $H=\re$ and
$A$ is the identity.

Also subproblem~(P2) is well understood, because on a \emph{fixed bounded}
time interval the degeneracy of the equation plays no role.
Local-in-time error estimates were proved by B.\ F.\ Esham and R.\
J.\ Weinacht in~\cite{ew}, then by the second author
in~\cite{k-cattaneo}, and finally by the authors
in~\cite[Appendix~A]{gg:k-PS} with optimal assumptions on initial
data.  The typical result is that $|A^{1/2}\roep(t)|\leq C\ep$ when
$(u_{0},u_{1})\in D(A^{3/2})\times D(A^{1/2})$, and we know that this
space is optimal if we look for estimates on $|A^{1/2}\roep(t)|$ of order
$\ep$, even in the linear case (see~\cite{gg:l-cattaneo}).

Subproblem~(P3) was solved by K.\ Nishihara and Y.\ Yamada~\cite{ny}.
They proved that (\ref{pbm:h-eq}), (\ref{pbm:h-data}) has a unique
global solution provided that $(u_{0},u_{1})\in\da\times\dau$ satisfy
the nondegeneracy assumption (\ref{hp:mdg}) and $\ep$ is small enough.
It is not known whether the smallness of $\ep$ is a necessary
condition. This remains the main open problem in the theory of
Kirchhoff equations, both dissipative and non-dissipative, both
degenerate and non-degenerate.

Subproblem~(P4) was first addressed in~\cite{ny}.  More recently, the
authors in~\cite{gg:k-decay} and~\cite{ghisi:decay} provided
\emph{optimal} decay estimates, showing that solutions of
(\ref{pbm:h-eq}), (\ref{pbm:h-data}) decay with the same rate of
solutions of the corresponding parabolic problem (see
also~\cite{mizu-ade,mizu-nc,ono-kyushu,ono-aa} for the case
$\gamma=1$).  The results have been recently extended
in~\cite{gg:w-dg} to equations with weak dissipation, namely with a
dissipative term of the form $b(t)\uep'(t)$, where $b(t)\to 0$ as
$t\to +\infty$.

Subproblem~(P5) was considered by the authors in~\cite{gg:k-PS},
with non-optimal convergence rates, and finally by the first
author~\cite{ghisi:error} with optimal convergence rates.

For the convenience of the reader, in section~\ref{sec:hystory} we
state all previous results needed in the sequel.

In this paper we concentrate on subproblem~(P6), namely on decay-error
estimates.  Estimates of this type were proved by R.\ Chill and A.\
Haraux~\cite{ch} in the case of \emph{linear equations}, and then by
H.\ Hashimoto and T.\ Yamazaki~\cite{yamazaki} for
\emph{nondegenerate} Kirchhoff equations.  Those results were
successively extended by T.\ Yamazaki~\cite{yamazaki-wd,yamazaki-cwd}
and by the authors~\cite{gg:w-ndg} to nondegenerate Kirchhoff
equations with weak dissipation.  The non-degenerate character of the
equation (namely strict hyperbolicity) seems to be \emph{essential} in
all previous approaches, which fail when applied to degenerate
equations.  This is the technical reason why subproblem~(P6) resisted
so far as an open problem.

In this paper we begin by showing that there is a deeper reason.
Indeed we show in Example~\ref{ex:no-de} that, without further
assumptions on initial data, the expected decay-error estimates are
actually false, even in the simple case where $H$ is a two dimensional
vector space.  By ``expected'' we mean decay-error estimates such as
(\ref{est:de-model}), where the decay-rate $\sigma(t)$ is the same as
in subproblem~(P4), and the convergence rate $\omega(\ep)$ is the same
as in subproblem~(P2) or subproblem~(P5).  The rigorous verification
of the counterexamples strongly relies on the asymptotic limits which
have been recently found in~\cite{ghisi:decay}.

Roughly speaking, the expected decay-error estimates are false
whenever the initial condition $u_{1}$ has a nonzero Fourier component
with respect to a frequency which is less than all frequencies
corresponding to nonzero components of $u_{0}$.  This motivates the
introduction of a special class of initial data where this cannot
happen (see Definition~\ref{hp:data}).  In Remark~\ref{rmk:data} we
show that this requirement on initial data is easily satisfied in many
concrete cases.

The main result of this paper is that in this class of initial data we
do have decay-error estimates for the degenerate problem.  Apart from
the special assumption, the regularity we require on initial data is
optimal, because it is the same which was optimal in the linear
nondegenerate case.  The convergence rates $\omega(\ep)$ are optimal,
because they are the same which appear in the local-in-time error
estimates of subproblem~(P2), or in the global-in-time error estimates
of subproblem~(P5).  The real surprise lies in the decay rate
$\sigma(t)$.  Indeed it turns out that $\roep(t)$ and $\rep(t)$ decay
\emph{faster} than $\uep(t)$ and $u(t)$ alone.

An improvement of decay rates has been observed also in~\cite{ch}
and~\cite{yamazaki,yamazaki-wd}, but in those cases it seems to
originate from different reasons.  Indeed in those examples it is
essential that the operator is not coercive, while in our case we have
improvement even if the operator is coercive.  Roughly speaking, our
improvement comes from the fact that our equation is in the same time
degenerate and nonlinear. In section~\ref{sec:heuristics} below we show a simple
toy model, based on ordinary differential equations of order one,
which gives a flavor of this aspect. The main point, both for the 
improvement and for the impossibility of expected decay-error 
estimates for general data, is that solutions of (\ref{eqn:ODE-model}) 
decay as $C(1+t)^{-1/(2\gamma)}$, where the constant $C$ depends on 
$\gamma$, but is independent of the initial condition.

Our result requires a new approach in order to take advantage of the
special assumptions on initial data.  The main idea is that in the
nonlinear degenerate case Fourier components corresponding to higher
frequencies decay faster.  As a consequence, in the limit as $t\to
+\infty$ the nonlinear terms $|A^{1/2}\uep(t)|^{2\gamma}$ or
$|A^{1/2}u(t)|^{2\gamma}$ do depend on the lowest frequency only.
This suggests to separate components corresponding to the lowest
frequency, and this is exactly what we do in Lemma~\ref{lemma:p-comp}
and then in section~\ref{sec:roep}, where we prove our basic
decay-error estimate on $\roep(t)$.  This is the nonlinear core of the
paper.

After the estimate on $\roep(t)$ has been established, the proof
becomes more standard.  We forget about components, and we regard both
(\ref{pbm:h-eq}) and (\ref{pbm:p-eq}) as linear equations where we
have frozen the nonlinear terms.  At this point we introduce weighted
versions of classical energies and we deduce all remaining
integral and pointwise estimates on $\roep(t)$, $\rep(t)$, and their
derivatives.

This paper is organized as follows. In section~\ref{sec:statements} we 
recall previous works, we state our main result, and we present some 
heuristics based on a toy model. In section~\ref{sec:proofs} we prove 
our main result. In section~\ref{sec:open} we state some open problems.

\setcounter{equation}{0}
\section{Statements}\label{sec:statements}

\subsection{Previous results}\label{sec:hystory}

In this section we recall some previous results needed in the sequel,
adapting them to the special nonlinear term which appears in
(\ref{pbm:h-eq}) and (\ref{pbm:p-eq}).

The first one answers what we called subproblem~(P1) in the
introduction.  It can be easily deduced from the theory developed
in~\cite{k-par}.  We recall that an operator $A$ is \emph{coercive} if
$$\inf\left\{\langle Au,u\rangle:u\in D(A),\ |u|=1\right\}>0.$$

\begin{thmbibl}[Parabolic problem: global existence and decay 
	estimates] 	\label{thm:hyst-p} 
	\hskip 0em plus 2em
	Let $H$ be a Hilbert space, and let $A$ be a nonnegative
	self-adjoint (unbounded) operator on $H$ with dense domain. Let 
	$\gamma\geq 1$ be a real number, and let $u_{0}\in D(A)$.
	
	Then we have the following conclusions.
	\begin{enumerate}
		\renewcommand{\labelenumi}{(\arabic{enumi})} 
		\item \emph{(Existence and uniqueness)}  Problem
		(\ref{pbm:p-eq}), (\ref{pbm:p-data}) has a unique global solution 
		$$u\in C^{1}\left([0,+\infty);H\right)\cap
		C^{0}\left([0,+\infty);\da\right).$$
	
		\item \emph{(Further regularity)}  If in addition $u_{0}$
		satisfies the nondegeneracy assumption (\ref{hp:mdg}), then
		the solution is non-stationary, and $u\in
		C^{\infty}\left((0,+\infty);D(A^{\alpha})\right)$ for every
		$\alpha\geq 0$.
	
		\item \emph{(Decay estimates)} Let us assume that the operator
		$A$ is coercive, that $u_{0}$ satisfies the nondegeneracy
		assumption (\ref{hp:mdg}), and that $u_{0}\in D(A^{k/2})$ for
		some integer $k\geq 2$.
		
		Then there exist positive constants $C_{1}$ and $C_{2}$ such
		that, for every positive integer $j\leq k$, we have that
		\begin{equation}
			\frac{C_{1}}{(1+t)^{1/\gamma}}\leq |A^{j/2}u(t)|^{2}\leq
			\frac{C_{2}}{(1+t)^{1/\gamma}} \quad\quad\forall t\geq 0.
			\label{th:dec-par}
		\end{equation}
	\end{enumerate}
\end{thmbibl}

The second result concerns subproblems~(P3) and~(P4). Existence and 
uniqueness were proved in~\cite{ny} (see 
also~\cite{gg:k-dissipative}), while decay estimates were proved in 
this form in~\cite{gg:k-decay}.

\begin{thmbibl}[Hyperbolic problem: global existence and decay estimates]
	\label{thm:hyst-h} 
	\hskip 0em plus 2em \mbox{} Let $H$ be a Hilbert space, and let
	$A$ be a nonnegative self-adjoint (unbounded) operator on $H$ with
	dense domain.  Let $\gamma\geq 1$ be a real number, and let us
	assume that the initial condition $(u_{0},u_{1})\in\da\times\dau$
	satisfies the nondegeneracy assumption (\ref{hp:mdg}).
	
	Then there exists $\ep_{0}>0$ for which the following conclusions
	hold true.
	\begin{enumerate}
		\renewcommand{\labelenumi}{(\arabic{enumi})} 
		
		\item \emph{(Existence and uniqueness)} For every
		$\ep\in(0,\ep_{0})$ we have that problem (\ref{pbm:h-eq}),
		(\ref{pbm:h-data}) has a unique global solution $\uep$ in the
		space
		\begin{equation}
			C^{2}\left([0,+\infty);H\right)\cap
			C^{1}\left([0,+\infty);\dau\right)\cap
			C^{0}\left([0,+\infty);\da\right).
			\label{reg:space}
		\end{equation}
	
		\item \emph{(Decay estimates)}  Let us assume in addition that
		the operator $A$ is coercive.  Then there exist positive
		constants $C_{1}$ and $C_{2}$ such that, for every
		$\ep\in(0,\ep_{0})$, we have that
		\begin{equation}
			\frac{C_{1}}{(1+t)^{1/\gamma}}\leq \auq{\uep(t)}\leq
			\frac{C_{2}}{(1+t)^{1/\gamma}} \quad\quad\forall t\geq 0,
			\label{est:auq}
		\end{equation}
		\begin{equation}
			\frac{C_{1}}{(1+t)^{1/\gamma}}\leq
			|A\uep(t)|^{2}\leq
			\frac{C_{2}}{(1+t)^{1/\gamma}}
			\quad\quad\forall t\geq 0,
			\label{est:au}
		\end{equation}
		\begin{equation}
			|\uep'(t)|^{2}\leq
			\frac{C_{2}}{(1+t)^{2+1/\gamma}}
			\quad\quad\forall t\geq 0.
			\label{est:u'}
		\end{equation}
	\end{enumerate}
\end{thmbibl}

The third and last result answers subproblem~(P5).  It follows from a
more general result proved in~\cite{ghisi:error} (see
also~\cite{gg:k-PS}), where also weak dissipation terms are
considered.

\begin{thmbibl}[Singular perturbation: global-in-time error estimates]
	\label{thm:hyst-sp}
	\hskip 0em plus 1em\mbox{}
	Let $H$, $A$, $\gamma$, $(u_{0},u_{1})$, $\ep_{0}$ be as in
	Theorem~\ref{thm:hyst-h}, and let $\uep(t)$, $u(t)$, $\tetep(t)$,
	$\roep(t)$, $\rep(t)$ be defined as usual.
	
	Let us assume that the operator $A$ is coercive.
	
	Then the following conclusions hold true.
	\begin{enumerate}
		\renewcommand{\labelenumi}{(\arabic{enumi})}
		\item  If in addition we assume that 
		$(u_{0},u_{1})\in\dat\times\dau$, then there exist
		$\ep_{1}\in(0,\ep_{0})$ and a constant $C$ such that for
		every $\ep\in(0,\ep_{1})$ we have that
		$$|\roep(t)|^{2}+|A^{1/2}\roep(t)|^{2}+ \ep(1+t)|\rep'(t)|^{2}
		\leq C\ep^{2}
		\quad\quad
		\forall t\geq 0,$$
		$$\int_{0}^{+\infty}\left(
		(1+t)|\rep'(t)|^{2}+
		\frac{|A^{1/2}\roep(t)|^{2}}{1+t}\right)dt
		\leq C\ep^{2}.$$
		
		\item  If in addition we assume that 
		$(u_{0},u_{1})\in D(A^{2})\times\da$, then there exist
		$\ep_{1}\in(0,\ep_{0})$ and a constant $C$ such that for
		every $\ep\in(0,\ep_{1})$ we have that
		$$|A\roep(t)|^{2}+(1+t)^{2}|\rep'(t)|^{2}
		\leq C\ep^{2}
		\quad\quad
		\forall t\geq 0,$$
		$$\int_{0}^{+\infty} \left((1+t)
		|A^{1/2}\rep'(t)|^{2}+
		\frac{|A\roep(t)|^{2}}{1+t}\right)dt\leq C\ep^{2}.$$
	\end{enumerate}
\end{thmbibl}

\begin{rmk}
	\begin{em}
		Decay-error estimates can be obtained by combining 
		Theorems~\ref{thm:hyst-p}, \ref{thm:hyst-h}, and 
		\ref{thm:hyst-sp} with standard inequalities. For example, 
		from previous results we know that
		$$|A^{1/2}\roep(t)|^{2}\leq 2\left(
		|A^{1/2}\uep(t)|^{2}+|A^{1/2}u(t)|^{2}\right)\leq
		\frac{K_{1}}{(1+t)^{1/\gamma}}
		\quad\mbox{and}\quad
		|A^{1/2}\roep(t)|^{2}\leq K_{2}\ep^{2}.$$
		
		Since $\min\{x,y\}\leq x^{\theta}y^{1-\theta}$ for every 
		$x\geq 0$, $y\geq 0$, $\theta\in(0,1)$, we have also that
		\begin{equation}
			|A^{1/2}\roep(t)|^{2}\leq\min\left\{
			\frac{K_{1}}{(1+t)^{1/\gamma}},K_{2}\ep^{2}\right\}\leq
			K_{3}\frac{\ep^{2(1-\theta)}}{(1+t)^{\theta/\gamma}}
			\quad\quad
			\forall t\geq 0.
			\label{est:de-trivial}
		\end{equation}
		
		These estimates are in general nonoptimal, both for the decay 
		rate, and for the convergence rate.
	\end{em}
\end{rmk}

\subsection{Notation and main result}\label{sec:main}

Taking into account the decay rates of Theorem~\ref{thm:hyst-h}, and 
the convergence rates of Theorem~\ref{thm:hyst-sp}, it is reasonable 
to expect decay-error estimates such as
\begin{equation}
	|A^{1/2}\roep(t)|^{2}\leq
	K\frac{\ep^{2}}{(1+t)^{1/\gamma}}
	\quad\quad
	\forall t\geq 0.
	\label{est:de-expected}
\end{equation}

The following example shows that such an estimate cannot be true for 
all initial data, even in finite dimension.

\begin{ex}\label{ex:no-de}
	\begin{em}
		Let us take $H:=\re^{2}$, and an operator $A$ with two 
		eigenvalues $\lambda_{0}^{2}<\lambda_{1}^{2}$, with 
		corresponding eigenvectors $e_{0}$ and $e_{1}$. Let us 
		consider the solutions of (\ref{pbm:h-eq}) and 
		(\ref{pbm:p-eq}) with initial data $u_{0}:=e_{1}$ and 
		$u_{1}:=e_{0}$. Let us write in components 
		$\uep(t)=u_{\ep,0}(t)e_{0}+u_{\ep,1}(t)e_{1}$, and 
		$u(t)=u_{0}(t)e_{0}+u_{1}(t)e_{1}$. Then it is easy to see 
		that $u_{0}(t)\equiv 0$. Moreover, from Theorem~3.3 
		of~\cite{ghisi:decay} we have that
		$$\lim_{t\to +\infty}(1+t)^{1/\gamma}|u_{\ep,0}(t)|^{2}=
		\frac{1}{\lambda_{0}^{2}}
		\frac{1}{(2\gamma\lambda_{0}^{2})^{1/\gamma}},$$
		and in particular the limit is different from 0 and 
		$\ep$-independent. It follows that
		$$\liminf_{t\to +\infty}(1+t)^{1/\gamma}|A^{1/2}\roep(t)|^{2}
		\geq \liminf_{t\to +\infty}(1+t)^{1/\gamma}
		\lambda_{0}^{2}|u_{\ep,0}(t)|^{2}=
		\frac{1}{(2\gamma\lambda_{0}^{2})^{1/\gamma}},$$
		which contradicts (\ref{est:de-expected}).
	\end{em}
\end{ex}

This example motivates the introduction of a class of initial data
where components of $u_{1}$ correspond to frequencies greater than or
equal to frequencies of components of $u_{0}$.  In order to state the
condition in a general form, we need some basic facts from the
spectral theory of operators, which we recall following~\cite{rudin}.

Let $E$ be the resolution of the identity associated with the 
operator $A$. For every measurable subset 
$J\subseteq[0,+\infty)$ we consider the space 
$H_{J}:=\mathcal{R}(E(J))$, namely the range of the projection 
operator $E(J)$, which is a closed subspace of $H$. For every 
$\mu>0$, we can therefore write $H$ as a direct sum
\begin{equation}
	H=H_{[0,\mu)}\oplus H_{[\mu,+\infty)}.
	\label{oplus}
\end{equation}

As a consequence, every vector $v\in H$ can be written in a unique 
way in the form $v=v_{\ell,\mu}+v_{h,\mu}$, with $v_{\ell,\mu}\in 
H_{[0,\mu)}$ and $v_{h,\mu}\in H_{[\mu,+\infty)}$. Here subscripts 
refer to low and high frequencies with respect to $\mu$. We also 
point out that
\begin{equation}
	\langle Av,v\rangle\geq \mu|v|^{2}
	\quad\quad
	\forall v\in\da\cap H_{[\mu,+\infty)}.
	\label{coercivity}
\end{equation}

In the case where $H$ admits a (finite or countable) orthonormal
system $\{e_{k}\}$ made by eigenvalues of $A$, and
$\{\lambda_{k}^{2}\}$ is the sequence of corresponding eigenvalues,
then $H_{J}$ is just the set of all $v\in H$ such that $\langle
v,e_{k}\rangle=0$ for every $k\in\n$ such that $\lambda_{k}^{2}\not\in
J$. Moreover in this case we have that
$$v_{\ell,\mu}:=\sum_{k:\lambda_{k}^{2}<\mu}\langle v,e_{k}\rangle
e_{k},
\hspace{3em}
v_{h,\mu}:=\sum_{k:\lambda_{k}^{2}\geq\mu}\langle v,e_{k}\rangle 
e_{k}.$$

We are now ready to introduce the class of initial data which is 
crucial for our decay-error estimates.

\begin{defn}[Assumption on initial data]\label{hp:data}
	\begin{em}
		Let $\nu>0$ and $\delta_{0}>1$ be two real numbers.  We say
		that a pair of initial conditions
		$(u_{0},u_{1})\in\da\times\dau$ satisfies the
		$(\nu,\delta_{0})$-assumption if 
		\begin{itemize}
			\item $\nu^{2}$ is an eigenvalue of $A$,
			
			\item the low frequency component
			$u_{0,\ell,\delta_{0}\nu^{2}}$ of $u_{0}$ is an
			eigenvector of $A$ (hence different from zero)
			corresponding to the eigenvalue $\nu^{2}$,
		
			\item  the low frequency component
			$u_{1,\ell,\delta_{0}\nu^{2}}$ of $u_{1}$ is a multiple 
			(possibly equal to zero) of the corresponding component
			of $u_{0}$, namely $u_{1,\ell,\delta_{0}\nu^{2}}=\beta 
			u_{0,\ell,\delta_{0}\nu^{2}}$ for some $\beta\in\re$.
		\end{itemize}
		
		In other words, $u_{0}$ is the sum of an eigenvector relative
		to $\nu^{2}$ and other components corresponding to frequencies
		greater than or equal to $\delta_{0}\nu^{2}$, while $u_{1}$ is
		the sum of a multiple (possibly equal to zero) of the same
		eigenvector and other components corresponding to frequencies
		greater than or equal to $\delta_{0}\nu^{2}$.
	\end{em}
\end{defn}

We point out that the $(\nu,\delta_{0})$-assumption implies (\ref{hp:mdg}).

\begin{rmk}\label{rmk:data}
	\begin{em}
		Let us assume that $H$ admits a (finite or countable)
		orthonormal system $\{e_{k}\}$ made by eigenvalues of $A$,
		relative to an increasing sequence of positive eigenvalues.
		Let $0<\lambda_{0}^{2}<\lambda_{1}^{2}$ be the two smallest
		eigenvalues.  Let us assume that $\lambda_{0}^{2}$ is simple,
		and let $e_{0}$ be a corresponding eigenvector.
		
		We point out that this assumption is always satisfied in the
		concrete case where $\Omega\subseteq\re^{n}$ is a connected
		bounded open set, $H:=L^{2}(\Omega)$, and $Au=-\Delta u$ with
		Dirichlet boundary conditions.  The interested reader is
		referred to Theorem~8.38 of~\cite{gt}.
		
		Let $(u_{0},u_{1})\in\da\times\dau$ be any initial condition
		such that $\langle u_{0},e_{0}\rangle\neq 0$.
		
		Then $(u_{0},u_{1})$ satisfies the 
		$(\nu,\delta_{0})$-assumption with $\nu:=\lambda_{0}$ and 
		$\delta_{0}:=\lambda_{1}^{2}/\lambda_{0}^{2}$.
	\end{em}
\end{rmk}

We can now state the main result of this paper.

\begin{thm}[Singular perturbation: decay-error estimates]
	\label{thm:main}
	Let $H$ be a Hil\-bert space, and let $A$ be a nonnegative
	self-adjoint (unbounded) operator on $H$ with dense domain.  Let
	$\gamma\geq 1$, $\nu>0$, and $\delta_{0}>1$ be real numbers.  Let
	$(u_{0},u_{1})\in\da\times\dau$ be a pair of initial conditions
	satisfying the $(\nu,\delta_{0})$-assumption.  Let $\ep_{0}$ be as
	in Theorem~\ref{thm:hyst-h}, and let $\uep(t)$, $u(t)$, $\tetep(t)$,
	$\roep(t)$, $\rep(t)$ be defined as usual.
	
	Let us set
	\begin{equation}
		\delta:=\min\{\delta_{0},2\gamma+1\},
		\label{defn:delta}
	\end{equation}
	and let us consider the function $\lambda:[0,+\infty)\to\re$ defined 
	by
	$$\lambda(t):=\left\{
		\begin{array}{ll}
			1 & \mbox{if }\delta<2\gamma+1,  \\
			\log(e+t) & \mbox{if }\delta=2\gamma+1.
		\end{array}
		\right.$$
	
	Then the following
	conclusions hold true.
	\begin{enumerate}
		\renewcommand{\labelenumi}{(\arabic{enumi})}
		\item If in addition we assume that
		$(u_{0},u_{1})\in\dat\times\dau$, then there exist
		$\ep_{1}\in(0,\ep_{0})$ and a constant $C$ such that for
		every $\ep\in(0,\ep_{1})$ we have that
		$$|\roep(t)|^{2}+|A^{1/2}\roep(t)|^{2}+ \ep(1+t)|\rep'(t)|^{2}
		\leq C\ep^{2}\frac{\lambda^{2}(t)}{(1+t)^{\delta/\gamma}}
		\quad\quad
		\forall t\geq 0,$$
		$$\int_{0}^{t}(1+s)^{2\delta/\gamma}\left((1+s)|\rep'(s)|^{2}+
		\frac{|A^{1/2}\roep(s)|^{2}}{1+s}\right)\,ds\leq 
		C\ep^{2}(1+t)^{\delta/\gamma}\lambda^{2}(t)
		\quad\quad
		\forall t\geq 0.$$
		
		\item If in addition we assume that $(u_{0},u_{1})\in
		D(A^{2})\times\da$, then there exist $\ep_{1}\in(0,\ep_{0})$
		and a constant $C$ such that for every $\ep\in(0,\ep_{1})$ we
		have that
		$$|A\roep(t)|^{2}+(1+t)^{2}|\rep'(t)|^{2}
		\leq C\ep^{2}\frac{\lambda^{2}(t)}{(1+t)^{\delta/\gamma}}
		\quad\quad
		\forall t\geq 0,$$
		$$\int_{0}^{t}(1+s)^{2\delta/\gamma}\left((1+s)|A^{1/2}\rep'(s)|^{2}+
		\frac{|A\roep(s)|^{2}}{1+s}\right)\,ds\leq 
		C\ep^{2}(1+t)^{\delta/\gamma}\lambda^{2}(t)
		\quad\quad
		\forall t\geq 0.$$
	\end{enumerate}
\end{thm}

\begin{rmk}\label{rmk:necessary}
	\begin{em}
		It is possible to show that there is no improvement of decay
		rates when $u_{1,\ell,\delta_{0}\nu^{2}}$ is not a multiple of
		$u_{0,\ell,\delta_{0}\nu^{2}}$.  The example is similar to
		Example~\ref{ex:no-de}, just with $\lambda_{0}=\lambda_{1}$.
		In this case a step of the proof of Theorem~3.3
		of~\cite{ghisi:decay} implies that $$\lim_{t\to
		+\infty}(1+t)^{1/\gamma}
		|u_{\ep,0}(t)|^{2}\neq 0.$$
		
		Here the limit could be $\ep$-dependent, but in any
		case this prevents $|\roep(t)|^{2}$ from decaying faster that
		$(1+t)^{1/\gamma}$.
		
		This shows that also the third condition in
		Definition~\ref{hp:data} is needed in order to have an
		improvement of decay rates.
	\end{em}
\end{rmk}

\subsection{Heuristics}\label{sec:heuristics}

A toy model for the singular perturbation problem is considering the 
difference between two solutions of the first order problem with two 
different initial data. The analogy is reasonable if we accept that the 
second order equation~(\ref{pbm:h-eq}) behaves as the first order 
equation~(\ref{pbm:p-eq}) when $\ep$ is small enough. Then we  
further simplify the model by taking $H:=\re$ and 
$A=\mbox{identity}$. Thus we have reduced ourselves to considering 
the difference between two different solutions of a first order 
ODE.

Despite of the dramatic simplification, the toy model still reveals a 
rich behavior. Indeed let us consider the following four examples.
\begin{enumerate}
	\renewcommand{\labelenumi}{(E\arabic{enumi})}
	\item  Let us examine equation $u'+u=0$ (linear and 
	nondegenerate). All solutions decay exponentially, and the 
	difference between two different solutions has the same decay rate of 
	the two solutions alone.

	\item Let us examine equation $u'+k(1+t)^{-1}u=0$ (linear and
	degenerate).  All solutions decay with a polynomial rate, and the
	difference between two different solutions decays with the same
	polynomial rate.

	\item Let us examine equation $u'+(1+|u|^{2\gamma})u=0$ (nonlinear and
	nondegenerate).  Once again solutions and differences between
	different solutions decay with the same (exponential) rate.

	\item Let us examine equation $u'+|u|^{2\gamma}u=0$ (nonlinear and
	degenerate).  Now solutions decay as $(1+t)^{-1/(2\gamma)}$, which
	is consistent with the decay rates in Theorem~\ref{thm:hyst-p} and
	Theorem~\ref{thm:hyst-h}.  On the contrary, the difference between
	two solutions with positive data decays as
	$(1+t)^{-(2\gamma+1)/(2\gamma)}$.  In other words, the decay rate
	of the difference is faster by a factor $(2\gamma+1)$.
\end{enumerate}

These examples seem to suggest that the improvement of decay rates
depends both on the nonlinear character and on the degeneracy of the
equation.  Last example suggests also that the factor $(2\gamma+1)$ in
the right-hand side of~(\ref{defn:delta}) is optimal.

Let us consider now the interaction between different Fourier
components.  For the sake of simplicity we take $H$, $A$,
$\lambda_{0}^{2}$, $\lambda_{1}^{2}$, $e_{0}$, and $e_{1}$ as in
Example~\ref{ex:no-de}.  Then we take the solution $u(t)$
of~(\ref{pbm:p-eq}) with initial condition $u(0)=e_{0}$, and the
solution $v(t)$ of~(\ref{pbm:p-eq}) with initial condition
$v(0)=e_{0}+e_{1}$.

It is easy to see that $u(t)$ has a unique component $u_{0}(t)e_{0}$,
whose coefficient satisfies
$u_{0}'(t)+\lambda_{0}^{2\gamma+2}u_{0}^{2\gamma+1}(t)=0$.  Once again
the solution decays as $(1+t)^{1/(2\gamma)}$, and an easy computation
shows that the nonlinear term is
\begin{equation}
	|A^{1/2}u(t)|^{2\gamma}=
	\lambda_{0}^{2\gamma}u_{0}^{2\gamma}(t)
	\sim\frac{1}{2\gamma\lambda_{0}^{2}}\frac{1}{(1+t)}.
	\label{heuristic}
\end{equation}

Now let us estimate $v(t)$.  It can be written in the form
$v(t)=v_{0}(t)e_{0}+v_{1}(t)e_{1}$, where $v_{0}(t)$ and $v_{1}(t)$
satisfy the system $v_{i}'(t)+\lambda_{i}^{2}c(t)v_{i}(t)=0$ (with
$i=0,1$), where
$$c(t)=|A^{1/2}v(t)|^{2\gamma}=
\left[\lambda_{0}^{2}v_{0}^{2}(t)+
\lambda_{1}^{2}v_{1}^{2}(t)\right]^{\gamma}.$$ 

We know from Theorem~\ref{thm:hyst-p} that $c(t)\sim(1+t)^{-1}$, hence
$v_{0}(t)$ and $v_{1}(t)$ decay with a polynomial rate with exponents
depending on $\lambda_{0}^{2}$ and $\lambda_{1}^{2}$.  In particular,
$v_{1}(t)$ decays faster than $v_{0}(t)$, so that in the limit it is
reasonable to assume that
$c(t)\sim\lambda_{0}^{2\gamma}v_{0}^{2\gamma}(t)$.

This ansatz uncouples the system, and therefore $v_{0}(t)$ becomes the 
solution of a single equation, the same solved by $u_{0}(t)$. This 
means that it is reasonable to assume that $u_{0}(t)\sim v_{0}(t)$, 
and $|A^{1/2}u(t)|^{2\gamma}\sim|A^{1/2}v(t)|^{2\gamma}$. At this 
point the difference $\rho(t)=v(t)-u(t)$ has a unique component 
$\rho_{1}(t)e_{1}=u_{1}(t)e_{1}$, so that 
$\rho_{1}'(t)+\lambda_{1}^{2}c(t)\rho_{1}(t)=0$. 

Setting $c(t)$ equal to the right-hand side of~(\ref{heuristic}), an 
easy computation shows that 
$$\rho_{1}(t)\sim
\frac{1}{(1+t)^{\lambda_{1}^{2}/(2\gamma\lambda_{0}^{2})}}
,$$ 
which means that there is an improvement of the decay rate equal to
$\lambda_{1}^{2}/\lambda_{0}^{2}$.  This suggests that the term
$\delta_{0}$ in~(\ref{defn:delta}) is optimal.

Our heuristic arguments are far from being proofs, even for the toy 
model of the difference between two solutions of the parabolic 
problem. Nevertheless, we hope that they can shed some light on the 
improvement of decay rates, and on the reason why it should depend on 
some $\delta$ defined by~(\ref{defn:delta}). 

\setcounter{equation}{0}
\section{Proofs}\label{sec:proofs}

This section is devoted to the proof of Theorem~\ref{thm:main}. In 
all proofs we set
\begin{equation}
	c(t):=|A^{1/2}u(t)|^{2\gamma},
	\hspace{3em}
	\cep(t):=|A^{1/2}\uep(t)|^{2\gamma}.
	\label{defn:c-cep}
\end{equation}

We also use that the corrector $\tetep(t)$, which is the solution of 
(\ref{pbm:tetep-eq}), (\ref{pbm:tetep-data}), is given by the explicit 
formula
\begin{equation}
	\tetep(t)=\ep w_{0}\left(1-e^{-t/\ep}\strut\right)
	\quad\quad
	\forall t\geq 0.
	\label{eqn:tetep}
\end{equation}

In many points we need to split vectors according to the orthogonal 
sum~(\ref{oplus}). In this case $v_{\ell,\mu}$ and $v_{h,\mu}$ denote 
the components of a certain vector $v\in H$, shortened to $v_{\ell}$ 
and $v_{h}$ when $\mu=\delta_{0}\nu^{2}$.

Due to our assumptions on initial data, all solutions lie in the 
space $H_{[\nu^{2},+\infty)}$. Therefore we can always assume, 
without loss of generality, that the operator is coercive, so that we 
can apply all the results stated in Theorems~\ref{thm:hyst-p},  
\ref{thm:hyst-h}, and~\ref{thm:hyst-sp}.

In all proofs, $k_{1}$, $k_{2}$, \ldots\ are real positive constants,
always independent of $\ep$ and $t$.  We restart the numeration of
constants in each proof.

\subsection{Preliminaries}

We recall some decay estimates for solutions of (\ref{pbm:h-eq}), 
(\ref{pbm:h-data}) which are needed in the sequel. The first one 
concerns the faster decay of components corresponding to high 
frequencies. A proof is contained in Theorem~3.1 and Theorem~3.3 
of~\cite{ghisi:decay}.

\begin{propbibl}[Faster decay for high frequencies]
	Let $H$, $A$, $\gamma$, $\nu$, $\delta_{0}$, $(u_{0},u_{1})$,
	$\ep_{0}$, $\uep(t)$ be as in Theorem~\ref{thm:main}.  Let $\mu>0$
	be a real number.
	
	Then there exist $\ep_{1}\in(0,\ep_{0})$, and a constant $M$ 
	(depending also on $\mu$), such that for every 
	$\ep\in(0,\ep_{1})$ we have that
	\begin{equation}
		|Au_{\ep,h,\mu}(t)|^{2}\leq
		\frac{M}{(1+t)^{\mu/(\nu^{2}\gamma)}}
		\quad\quad
		\forall t\geq 0,
		\label{th:mist-Auh}
	\end{equation}
	\begin{equation}
		|u_{\ep,h,\mu}'(t)|^{2}\leq
		\frac{M}{(1+t)^{2+\mu/(\nu^{2}\gamma)}}
		\quad\quad
		\forall t\geq 0.
		\label{th:mist-u'}
	\end{equation}
\end{propbibl}

The second result concerns the decay of second derivatives.  The
estimate deals with low frequencies, and it follows from Theorem~3.3
and Proposition~4.3 of~\cite{ghisi:decay}.  We point out that an
analogous estimate holds true without restricting to low frequencies
provided that initial data are more regular, namely $(u_{0},u_{1})\in
D(A^{2})\times\dat$, or with a constant $M$ which depends also
on~$\ep$.

\begin{propbibl}[Decay for low frequencies of second 
	derivatives]\label{prop:mist-u''}
	Let $H$, $A$, $\gamma$, $\nu$, $\delta_{0}$, $(u_{0},u_{1})$,
	$\ep_{0}$, $\uep(t)$, $\tetep(t)$ be as in Theorem~\ref{thm:main}.
	Let $\mu>0$ be a real number.
	
	Then there exist $\ep_{1}\in(0,\ep_{0})$, and a constant $M$ 
	(depending also on $\mu$), such that for every 
	$\ep\in(0,\ep_{1})$ we have that
	\begin{equation}
		|u_{\ep,\ell,\mu}''(t)-\theta_{\ep,\ell,\mu}''(t)|^{2}\leq
		\frac{M}{(1+t)^{4+1/\gamma}}
		\quad\quad
		\forall t\geq 0.
		\label{th:mist-u''}
	\end{equation}
\end{propbibl}

Now we state and prove two results for ordinary differential
equations.  The first one is a simple comparison principle, which has
already been used in similar forms
in~\cite{ghisi:error,ghisi:decay,gg:k-dissipative,gg:k-decay,gg:k-PS,gg:w-ndg,gg:w-dg}.

\begin{lemma}\label{lemma:ODE-solito}
	Let $\psi:[0,+\infty)\to(0,+\infty)$ be a nondecreasing function 
	of class $C^{1}$. Let $M$ be a positive constant, and let 
	$z:[0,+\infty)\to[0,+\infty)$ be a function of class $C^{1}$ such 
	that $z(0)=0$, and
	\begin{equation}
		z'(t)\leq-M\sqrt{z(t)}\left(\sqrt{z(t)}-\psi(t)\right)
		\quad\quad
		\forall t\geq 0.
		\label{hp:ODE}
	\end{equation}
	
	Then we have that $z(t)\leq\psi^{2}(t)$ for every $t\geq 0$.
\end{lemma}

\paragraph{\emph{\textmd{Proof}}}

Let us consider the differential equation
$y'=-M\sqrt{y}\left(\sqrt{y}-\psi(t)\right)$.
Assumption~(\ref{hp:ODE}) is equivalent to say that $z(t)$ is a
subsolution.  On the other hand, due to the monotonicity of $\psi(t)$,
it is easy to check that $w(t):=\psi^{2}(t)$ is a supersolution of the
same equation.  Since $z(0)=0<w(0)$, the conclusion follows from the
standard comparison principle between subsolutions and
supersolutions.\qed
\medskip

The second lemma is a comparison result for a more complex
differential inequality.  The assumptions on the coefficients are
exactly those which are satisfied in section~\ref{sec:roep}, where
this lemma plays a crucial role.

\begin{lemma}\label{lemma:main-ODE}
	Let $\ep_{0}>0$, let $\lambda:[0,+\infty)\to[1,+\infty)$ be a
	continuous function, and let
	$\psi_{i}:(0,\ep_{0})\times[0,+\infty)\to\re$ (with $i=1,2,3,4$)
	be continuous functions.  
	
	Let us assume that there exist constants
	$M_{1}$, \ldots, $M_{5}$ such that for every $\ep\in(0,\ep_{0})$
	we have that
	\begin{eqnarray}
		& \psi_{1}(\ep,t)\geq 0
		\quad\quad\forall t\geq 0, & 
	   \label{hp:psi0}  \\
	   \noalign{\vspace{1ex}}
	   & \log\lambda(t)\leq
	   \displaystyle{\int_{0}^{t}\psi_{1}(\ep,s)\,ds}\leq
	   M_{1}+\log\lambda(t)
	   \quad\quad\forall t\geq 0, & 
	  \label{hp:psi1}  \\
		\noalign{\vspace{1ex}}
		 & \displaystyle{\int_{0}^{+\infty}
		 |\psi_{2}(\ep,s)|\cdot\lambda^{3}(s)\,ds}\leq
		 M_{2}, & 
		\label{hp:psi2}  \\
		\noalign{\vspace{1ex}}
		 & \displaystyle{\int_{0}^{t}
		 \frac{|\psi_{3}(\ep,s)|}{\lambda(s)}\,ds}\leq
		 (M_{3}+M_{4}\lambda(t))\ep^{2}
		 \quad\quad\forall t\geq 0, & 
		\label{hp:psi3}  \\
		\noalign{\vspace{1ex}}
		 & \displaystyle{\left|\int_{0}^{t}
		 \psi_{4}(\ep,s)\,ds\right|}\leq
		 M_{5}\ep^{2}
		 \quad\quad\forall t\geq 0. & 
		\label{hp:psi4}  
	\end{eqnarray}
	
	For every $\ep\in(0,\ep_{0})$, let 
	$z_{\ep}:[0,+\infty)\to[0,+\infty)$ be a function of class 
	$C^{1}$ such that $z_{\ep}(0)=0$, and
	\begin{equation}
		z_{\ep}'(t)\leq\psi_{1}(\ep,t)z_{\ep}(t)+
		\psi_{2}(\ep,t)\left[z_{\ep}(t)\right]^{3/2}+
		\psi_{3}(\ep,t)+\psi_{4}(\ep,t)
		\quad\quad
		\forall t\geq 0.
		\label{eqn:zep}
	\end{equation}
	
	Then there exist $\ep_{1}\in(0,\ep_{0})$, and a constant $M_{6}$
	such that for every $\ep\in(0,\ep_{1})$ we have that
	\begin{equation}
		z_{\ep}(t)\leq M_{6}\ep^{2}\lambda^{2}(t)
		\quad\quad
		\forall t\geq 0.
		\label{th:lemma-main}
	\end{equation}
\end{lemma}

\paragraph{\emph{\textmd{Proof}}}

For every $i\in\{1,2,3,4\}$, and every $\ep\in(0,\ep_{0})$, let us set
$$\Psi_{i}(\ep,t):=\int_{0}^{t}\psi_{i}(\ep,s)\,ds
\quad\quad
\forall t\geq 0.$$

For the sake of simplicity, when no confusion is possible we omit the
dependence on $\ep$, and sometimes also the dependence on $t$, when
writing $z_{\ep}(t)$, $\psi_{i}(\ep,t)$, $\Psi_{i}(\ep,t)$.  In any
case all constants we introduce are independent of $\ep$ and $t$.

From differential inequality (\ref{eqn:zep}) we have that
\begin{eqnarray*}
	\left[e^{-\Psi_{1}}z\right]' & \leq &
	e^{-\Psi_{1}}\psi_{2}\,z^{3/2}+
	e^{-\Psi_{1}}\psi_{3}+e^{-\Psi_{1}}\psi_{4} \\
	 & = & e^{-\Psi_{1}}\psi_{2}\,z^{3/2}+
	e^{-\Psi_{1}}\psi_{3}+
	\left[e^{-\Psi_{1}}\Psi_{4}\right]'+
	e^{-\Psi_{1}}\Psi_{4}\psi_{1}.
\end{eqnarray*}

Integrating in $[0,t]$, and exploiting the initial condition 
$z(0)=0$, we obtain that
\begin{eqnarray}
	z(t) & \leq & e^{\Psi_{1}(t)}\int_{0}^{t}
	e^{-\Psi_{1}(s)}\psi_{2}(s)[z(s)]^{3/2}\,ds+
	e^{\Psi_{1}(t)}\int_{0}^{t}
	e^{-\Psi_{1}(s)}\psi_{3}(s)\,ds 
	\nonumber \\
	 &  & +\Psi_{4}(t)+e^{\Psi_{1}(t)}\int_{0}^{t}
	e^{-\Psi_{1}(s)}\psi_{1}(s)\Psi_{4}(s)\,ds 
	\nonumber  \\
	& =: & I_{1}+I_{2}+I_{3}+I_{4}. 
	\label{est:zep}
\end{eqnarray}

Let us estimate the four terms. From (\ref{hp:psi1}), and the fact that 
$\lambda(t)\geq 1$, we have that
\begin{equation}
	I_{1}\leq e^{M_{1}}\lambda(t)\int_{0}^{t}
	|\psi_{2}(s)|\cdot [z(s)]^{3/2}\,ds.
	\label{est:lemma-I1}
\end{equation}

Exploiting (\ref{hp:psi1}), (\ref{hp:psi3}), and the fact that 
$\lambda(t)\leq\lambda^{2}(t)$, we obtain that
\begin{equation}
	I_{2}\leq e^{M_{1}}\lambda(t)\int_{0}^{t}
	\frac{|\psi_{3}(s)|}{\lambda(s)}\,ds\leq
	e^{M_{1}}\lambda(t)\left(
	M_{3}+M_{4}\lambda(t)\right)\ep^{2}\leq
	k_{1}\ep^{2}\lambda^{2}(t).
	\label{est:lemma-I2}
\end{equation}

Moreover assumption (\ref{hp:psi4}) is equivalent to say that
\begin{equation}
	I_{3}\leq\left|\Psi_{4}(t)\right|\leq M_{5}\ep^{2}.
	\label{est:lemma-I3}
\end{equation}

Finally, from (\ref{hp:psi0}) and (\ref{hp:psi4}) we have that
\begin{equation}
	I_{4}\leq e^{M_{1}}\lambda(t)\int_{0}^{t}
	\left|\Psi_{4}(s)\right|\cdot\psi_{1}(s)e^{-\Psi_{1}(s)}\,ds
	\leq k_{2}\ep^{2}\lambda(t)
	\int_{0}^{t}\psi_{1}(s)e^{-\Psi_{1}(s)}\,ds\leq
	k_{2}\ep^{2}\lambda^{2}(t).
	\label{est:lemma-I4}
\end{equation}

Plugging (\ref{est:lemma-I1}) through (\ref{est:lemma-I4}) into 
(\ref{est:zep}) we obtain that
\begin{equation}
	z(t)\leq e^{M_{1}}\lambda(t)\int_{0}^{t}
	|\psi_{2}(s)|\cdot[z(s)]^{3/2}\,ds+
	k_{3}\ep^{2}\lambda^{2}(t)
	\quad\quad
	\forall t\geq 0.
	\label{est:lemma-main}
\end{equation}

Now let us choose $\ep_{1}$ small enough so that
$$M_{2}e^{M_{1}}3^{3/2}k_{3}^{1/2}\ep_{1}\leq 1,$$ 
and then let us set
$$T_{\ep}:=\sup\left\{t\geq 0:z_{\ep}(\tau)\leq
3k_{3}\ep^{2}\lambda^{2}(\tau)\quad\forall\tau\in[0,t]\right\}.$$

We claim that $T_{\ep}=+\infty$ for every $\ep\in(0,\ep_{1})$, which 
implies (\ref{th:lemma-main}). To this end, let us assume by 
contradiction that $T_{\ep}<+\infty$ for some $\ep\in(0,\ep_{1})$. 
Then we have that $T_{\ep}>0$ because $z_{\ep}(0)=0$, and
\begin{equation}
	z_{\ep}(t)\leq 3k_{3}\ep^{2}\lambda^{2}(t)
	\quad\quad
	\forall t\in[0,T_{\ep}],
	\label{eqn:zT2}
\end{equation}
\begin{equation}
	z_{\ep}(T_{\ep})=3k_{3}\ep^{2}\lambda^{2}(T_{\ep}).
	\label{eqn:zT1}
\end{equation}

Setting $t=T_{\ep}$ in (\ref{est:lemma-main}), and exploiting
(\ref{eqn:zT2}), we obtain that
$$z_{\ep}(T_{\ep})\leq e^{M_{1}}\lambda(T_{\ep})\int_{0}^{T_{\ep}}
|\psi_{2}(\ep,s)|\cdot 3^{3/2}k_{3}^{3/2}\ep^{3}\lambda^{3}(s)\,ds+
k_{3}\ep^{2}\lambda^{2}(T_{\ep}).$$

Exploiting (\ref{hp:psi2}), and inequalities
$\lambda(t)\leq\lambda^{2}(t)$ and $\ep\leq\ep_{1}$, we finally deduce that
\begin{eqnarray*}
	z_{\ep}(T_{\ep}) & \leq & e^{M_{1}}\lambda(T_{\ep})\cdot
	3^{3/2}k_{3}^{3/2}\ep^{3}\cdot M_{2} +
	k_{3}\ep^{2}\lambda^{2}(T_{\ep})   \\
	 & \leq & \left(M_{2}e^{M_{1}}
	 3^{3/2}k_{3}^{1/2}\ep_{1}+1 \right)k_{3}\ep^{2}\lambda^{2}(T_{\ep}) \\
	 & \leq & 2k_{3}\ep^{2}\lambda^{2}(T_{\ep}),
\end{eqnarray*}
which contradicts (\ref{eqn:zT1}).\qed

\subsection{Estimates on the parabolic equation}

In this section we collect the estimates on the parabolic equation, 
not contained in Theorem~\ref{thm:hyst-p}, which are needed in the 
proof of our main result.

The first one is an estimate on second derivatives. In particular, 
estimate~(\ref{th:par-v''-point}) is in some sense the parabolic 
counterpart of~(\ref{th:mist-u''}). Here we do not need to restrict 
to low frequencies because $u_{0}\in D(A^{2})$.

\begin{lemma}[Parabolic problem: estimates on second derivative]
	Let $H$, $A$, $\gamma$, $u_{0}$ be as in 
	Theorem~\ref{thm:hyst-p}, and let $u(t)$ be the corresponding 
	solution of problem~(\ref{pbm:p-eq}), (\ref{pbm:p-data}). 
	
	Then for every $\delta\in(0,2\gamma+1]$ we have the following
	conclusions.
	\begin{enumerate}
		\renewcommand{\labelenumi}{(\arabic{enumi})} 
		\item  If $u_{0}\in\dat$, then there exists a constant $M$ 
		such that
		\begin{equation}
			\int_{0}^{t}(1+s)^{1+2\delta/\gamma}|u''(s)|^{2}\,ds\leq
			M(1+t)^{\delta/\gamma}
			\quad\quad
			\forall t\geq 0.
			\label{th:par-v''}
		\end{equation}
	
		\item  If $u_{0}\in D(A^{2})$, then there exists a constant $M$ 
		such that
		\begin{equation}
			\int_{0}^{t}(1+s)^{1+2\delta/\gamma}|A^{1/2}u''(s)|^{2}\,ds\leq
			M(1+t)^{\delta/\gamma}
			\quad\quad
			\forall t\geq 0,
			\label{th:par-auq-v''}
		\end{equation}
		\begin{equation}
			|u''(t)|^{2}\leq\frac{M}{(1+t)^{4+1/\gamma}}
			\quad\quad
			\forall t\geq 0.
			\label{th:par-v''-point}
		\end{equation}
	\end{enumerate}
\end{lemma}

\paragraph{\textmd{\emph{Proof}}}

Let us set for simplicity $\eta:=\delta/\gamma$. Since $\delta\leq 2\gamma+1$, 
it follows that $\eta\leq 2+1/\gamma$, hence for every $t\geq 0$ we have 
that
\begin{equation}
	\int_{0}^{t}\frac{(1+s)^{2\dg-3}}{(1+s)^{1/\gamma}}\,ds=
	\int_{0}^{t}\frac{(1+s)^{2\dg+1}}{(1+s)^{4+1/\gamma}}\,ds\leq
	\int_{0}^{t}(1+s)^{\dg-1}\,ds\leq
	\frac{1}{\dg}(1+t)^{\dg}.
	\label{est:int-eta}
\end{equation}

\subparagraph{\textmd{\emph{Basic integral estimates}}}

We prove that when $u_{0}\in\dat$ we have that
\begin{equation}
	\int_{0}^{t}(1+s)^{2\dg-3}|A^{2}u(s)|^{2}\,ds\leq 
	k_{1}(1+t)^{\dg}
	\quad\quad
	\forall t\geq 0,
	\label{est:int-par-2}
\end{equation}
and when $u_{0}\in D(A^{2})$ we have that
\begin{equation}
	\int_{0}^{t}(1+s)^{2\dg-3}|A^{5/2}u(s)|^{2}\,ds\leq k_{2}(1+t)^{\dg}
	\quad\quad
	\forall t\geq 0.
	\label{est:int-par-5/2}
\end{equation}

To this end, an easy calculation shows that
$$\frac{d}{dt}
\left(\frac{1}{2}(1+t)^{2\dg-2}|A^{3/2}u(t)|^{2}\right)+
c(t)(1+t)^{2\dg-2}|A^{2}u(t)|^{2}=
(\dg-1)(1+t)^{2\dg-3}|A^{3/2}u(t)|^{2}.$$

Now we integrate in $[0,t]$, and then we apply the estimate from above
in (\ref{th:dec-par}) with $j=3$, and finally estimate
(\ref{est:int-eta}).  We obtain that
\begin{eqnarray*}
	\lefteqn{\hspace{-3em}\frac{1}{2}(1+t)^{2\dg-2}|A^{3/2}u(t)|^{2}+
	\int_{0}^{t}c(s)(1+s)^{2\dg-2}|A^{2}u(s)|^{2}\,ds} \\
	 & = & \frac{1}{2}|A^{3/2}u_{0}|^{2}+
	 (\dg-1)\int_{0}^{t}(1+s)^{2\dg-3}|A^{3/2}u(s)|^{2}\,ds \\
	 & \leq & \frac{1}{2}|A^{3/2}u_{0}|^{2}+
	 k_{3}\int_{0}^{t}(1+s)^{2\dg-3}(1+s)^{-1/\gamma}\,ds\\
	 & \leq & k_{4}(1+t)^{\dg}.
\end{eqnarray*}

Applying the estimate from below in (\ref{th:dec-par}) with $j=1$, 
we therefore deduce that
$$\int_{0}^{t}(1+s)^{2\dg-3}|A^{2}u(s)|^{2}\,ds \leq
k_{5}\int_{0}^{t}c(s)(1+s)^{2\dg-2}|A^{2}u(s)|^{2}\,ds
\leq k_{6}(1+t)^{\dg},$$
which proves (\ref{est:int-par-2}).  The proof of
(\ref{est:int-par-5/2}) is analogous (one just needs to add 1/2 to all
powers of the operator $A$).

\subparagraph{\textmd{\emph{Estimates on second derivatives}}}

Taking the time derivative of (\ref{pbm:p-eq}) we find 
that 
$$u''(t)=-c'(t)Au(t)-c(t)Au'(t)=
2\gamma|A^{1/2}u(t)|^{4\gamma-2}|Au(t)|^{2} Au(t)+
|A^{1/2}u(t)|^{4\gamma}A^{2}u(t)$$
for every $t>0$, hence
$$|u''(t)|^{2}\leq k_{7}|A^{1/2}u(t)|^{8\gamma-4}\cdot|Au(t)|^{6}+
k_{8}|A^{1/2}u(t)|^{8\gamma}\cdot|A^{2}u(t)|^{2},$$
$$|A^{1/2}u''(t)|^{2}\leq k_{7}|A^{1/2}u(t)|^{8\gamma-4}\cdot|Au(t)|^{4}
\cdot|A^{3/2}u(t)|^{2}+ k_{8}|A^{1/2}u(t)|^{8\gamma}
\cdot|A^{5/2}u(t)|^{2}.$$

Since in any case $u_{0}\in\dat$, we can apply (\ref{th:dec-par}) with $j=1,2,3$. 
We obtain that
\begin{equation}
	|u''(t)|^{2}\leq\frac{k_{9}}{(1+t)^{4+1/\gamma}}+
	\frac{k_{10}}{(1+t)^{4}}|A^{2}u(t)|^{2},
	\label{est:u''}
\end{equation}
\begin{equation}
	|A^{1/2}u''(t)|^{2}\leq\frac{k_{11}}{(1+t)^{4+1/\gamma}}+
	\frac{k_{12}}{(1+t)^{4}}|A^{5/2}u(t)|^{2}.
	\label{est:u''-auq}
\end{equation}

If $u_{0}\in\dat$, then (\ref{th:par-v''}) follows from
(\ref{est:u''}), (\ref{est:int-eta}), and (\ref{est:int-par-2}).  

If $u_{0}\in D(A^{2})$, then (\ref{th:par-auq-v''}) follows from
(\ref{est:u''-auq}), (\ref{est:int-eta}), and (\ref{est:int-par-5/2}).
Finally, (\ref{th:par-v''-point}) follows from (\ref{est:u''}) and
(\ref{th:dec-par}) with $j=4$.\qed
\medskip

In the second result we take a solution of the parabolic problem, and
we estimate its components with respect to low and high frequencies. 
In particular, estimate~(\ref{th:par-gain}) is the parabolic 
counterpart of~(\ref{th:mist-Auh}) in the special case 
$\mu=\delta_{0}\nu^{2}$.

We assume that the initial datum $u_{0}\in D(A)$ has the same
structure required in Definition~\ref{hp:data}.  This means that there
exist $\nu>0$, $\delta_{0}>1$, and a decomposition
$u_{0}=u_{0,\ell}+u_{0,h}$, where $\nu^{2}$ is an eigenvalue of $A$,
$u_{0,\ell}\neq 0$ is an eigenvector relative to $\nu^{2}$, and $u_{0,h}\in
H_{[\delta_{0}\nu^{2},+\infty)}$.

In this case the solution $u(t)$ of problem~(\ref{pbm:p-eq}),
(\ref{pbm:p-data}) can be written in the form
$u(t)=u_{\ell}(t)+u_{h}(t)$, where $u_{\ell}(t)$ and $u_{h}(t)$ are
the solutions of the linear problems
\begin{equation}
	u_{\ell}'(t)+c(t)Au_{\ell}(t)=0, \quad\quad
	u_{\ell}(0)=u_{0,\ell},
	\label{eqn:ul}
\end{equation}
\begin{equation}
	u_{h}'(t)+c(t)Au_{h}(t)=0,
	\quad\quad
	u_{h}(0)=u_{0,h},
	\label{eqn:uh}
\end{equation}
where of course $c(t)$ is given by (\ref{defn:c-cep}).

\begin{lemma}[Parabolic problem: estimates on low and high frequencies]\label{lemma:p-comp}
	\hskip 0em plus 2em
	Let $H$, $A$, $\gamma$ be as in Theorem~\ref{thm:hyst-p}. Let 
	$\nu$, $\delta_{0}$, $u_{0}=u_{0,\ell}+u_{0,h}$, and 
	$u(t)=u_{\ell}(t)+u_{h}(t)$ be as above. 
	Let us set
	\begin{equation}
		\Phi(t):=\nu^{2}\gamma\int_{0}^{t}|A^{1/2}u_{\ell}(s)|^{2\gamma}\,ds
		\quad\quad
		\forall t\geq 0.
		\label{defn:Phi}
	\end{equation}
	
	Then there exist positive constants $M_{1}$, $M_{2}$, $M_{3}$ 
	such that
	\begin{equation}
		|Au_{h}(t)|^{2}\leq\frac{M_{1}}{(1+t)^{\delta_{0}/\gamma}}
		\quad\quad
		\forall t\geq 0,
		\label{th:par-gain}
	\end{equation}
	\begin{equation}
		M_{2}(1+t)\leq e^{2\Phi(t)}\leq M_{3}(1+t)
		\quad\quad
		\forall t\geq 0.
		\label{th:par-Phi}
	\end{equation}
\end{lemma}

\paragraph{\textmd{\emph{Proof}}}

For every $t\geq 0$ let us set
$$\begin{array}{ccc}
	\displaystyle{C(t):=\int_{0}^{t}|A^{1/2}u(s)|^{2\gamma}\,ds,} & 
	\quad\quad &
	\displaystyle{C_{\ell}(t):=\int_{0}^{t}|A^{1/2}\ul(s)|^{2\gamma}\,ds,}  \\
	\noalign{\vspace{2ex}}
	y(t):=e^{2\nu^{2}\gamma C(t)}, & & 
	y_{\ell}(t):=e^{2\nu^{2}\gamma C_{\ell}(t)}=e^{2\Phi(t)}.
\end{array}$$

Since $u_{0,\ell}$ is an eigenvector of $A$, it is easy to see that 
the solution of (\ref{eqn:ul}) is given by the explicit formula
\begin{equation}
	\ul(t)=u_{0,\ell}e^{-\nu^{2}C(t)}
	\quad\quad
	\forall t\geq 0.
	\label{formula-ul}
\end{equation}

\subparagraph{\textmd{\emph{Estimate from below for $y(t)$}}}

We claim that
\begin{equation}
	y(t)\geq k_{1}(1+t)
	\quad\quad
	\forall t\geq 0.
	\label{est:y}
\end{equation}

Indeed from (\ref{formula-ul}) we have that 
$|A^{1/2}\ul(t)|^{2\gamma}=|A^{1/2}u_{0,\ell}|^{2\gamma}\cdot 
e^{-2\nu^{2}\gamma C(t)}$, hence
$$y'(t)=2\nu^{2}\gamma|A^{1/2}u(t)|^{2\gamma}\cdot e^{2\nu^{2}\gamma C(t)}
\geq 2\nu^{2}\gamma|A^{1/2}\ul(t)|^{2\gamma}\cdot e^{2\nu^{2}\gamma C(t)}
=2\nu^{2}\gamma|A^{1/2}u_{0,\ell}|^{2\gamma},$$
from which (\ref{est:y}) immediately follows.

\subparagraph{\textmd{\emph{Estimate on high frequencies}}}

Thanks to~(\ref{eqn:uh}) and (\ref{coercivity}) with 
$\mu=\delta_{0}\nu^{2}$, we have that
$$\frac{d}{dt}|A^{j/2}u_{h}|^{2}=
2\langle A^{j/2}u_{h},A^{j/2}u_{h}'\rangle=
-2c(t)\langle A^{j/2}u_{h},AA^{j/2}u_{h}\rangle\leq
-2\delta_{0}\nu^{2}c(t)|A^{j/2}u_{h}|^{2}$$
for every $j=1,2$, and every $t>0$. Integrating in $[0,t]$, and 
exploiting (\ref{est:y}), we obtain that
\begin{equation}
	|A^{j/2}u_{h}(t)|^{2}\leq|A^{j/2}u_{0,h}|^{2}e^{-2\delta_{0}\nu^{2}C(t)}
	\leq\frac{k_{2}}{(1+t)^{\delta_{0}/\gamma}}
	\quad\quad
	\forall t\geq 0.
	\label{est:par-hf}
\end{equation}

Estimate (\ref{est:par-hf}) with $j=2$ is exactly (\ref{th:par-gain}).

\subparagraph{\textmd{\emph{Estimate on $C(t)-C_{\ell}(t)$}}}

We claim that 
\begin{equation}
	0\leq C(t)-C_{\ell}(t)\leq k_{3}
	\quad\quad
	\forall t\geq 0.
	\label{est:Cell}
\end{equation}

The estimate from below is trivial. In order to prove the estimate 
from above, let us consider the well known inequality
$$0\leq (x+y)^{\gamma}-x^{\gamma}\leq\gamma(x+y)^{\gamma-1}y
\quad\quad
\forall x\geq 0,\ \forall y\geq 0.$$

Setting $x:=|A^{1/2}u_{\ell}(t)|^{2}$, and $y:=|A^{1/2}u_{h}(t)|^{2}$, we 
obtain that
\begin{eqnarray*}
	|A^{1/2}u(t)|^{2\gamma}-|A^{1/2}\ul(t)|^{2\gamma} & = & 
	\left(|A^{1/2}\ul(t)|^{2}+|A^{1/2}\uh(t)|^{2}\right)^{\gamma}-
	\left(|A^{1/2}\ul(t)|^{2}\right)^{\gamma}\\
	 & \leq & \gamma|A^{1/2}u(t)|^{2(\gamma-1)}|A^{1/2}u_{h}(t)|^{2}.
\end{eqnarray*}

Exploiting (\ref{th:dec-par}) with $j=1$, and (\ref{est:par-hf}) with 
$j=1$, we obtain that
$$|A^{1/2}u(t)|^{2\gamma}-|A^{1/2}\ul(t)|^{2\gamma}\leq 
\frac{k_{4}}{(1+t)^{1+(\delta_{0}-1)/\gamma}}.$$

Since $\delta_{0}>1$, integrating in $[0,t]$ we deduce the estimate 
from above in (\ref{est:Cell}).

\subparagraph{\textmd{\emph{Estimate on $y_{\ell}(t)$}}}

We prove that 
\begin{equation}
	k_{5}(1+t)\leq y_{\ell}(t)\leq k_{6}(1+t)
	\quad\quad
	\forall t\geq 0,
	\label{est:z}
\end{equation}
which is exactly (\ref{th:par-Phi}). Indeed we have that
$$y_{\ell}'(t)=2\nu^{2}\gamma|A^{1/2}\ul(t)|^{2\gamma}\cdot
e^{2\nu^{2}\gamma C_{\ell}(t)}=
2\nu^{2}\gamma|A^{1/2}u_{0,\ell}|^{2\gamma}\cdot
e^{2\nu^{2}\gamma \left(C_{\ell}(t)-C(t)\right)},$$
so that from (\ref{est:Cell}) we deduce that $0<k_{7}\leq y_{\ell}'(t)\leq 
k_{8}$ for every $t\geq 0$. 

Integrating in $[0,t]$ we obtain (\ref{est:z}).\qed

\subsection{Proof of key decay-error estimate}\label{sec:roep}

This section is the key step in the proof of Theorem~\ref{thm:main}. 
Here we show that
\begin{equation}
	|\roep(t)|^{2}\leq 
	k_{1}\ep^{2}\frac{\lambda^{2}(t)}{(1+t)^{\delta/\gamma}}
	\quad\quad
	\forall t\geq 0.
	\label{th:roep}
\end{equation}

Let $c(t)$, $\cep(t)$, and components of vectors be defined as in the
first paragraph of section~\ref{sec:proofs}.  Let
$\Phi:[0,+\infty)\to[0,+\infty)$ be the function defined by
(\ref{defn:Phi}).  Let us set for simplicity $\eta:=\delta/\gamma$,
and let
$$z_{\ep}(t):=\frac{1}{2}|\roep(t)|^{2}e^{2\eta\Phi(t)}.$$

We claim that $z_{\ep}(t)$ satisfies a differential inequality as in
Lemma~\ref{lemma:main-ODE}.  If we prove this claim, then from that
lemma it follows that
\begin{equation}
	z_{\ep}(t)\leq k_{2}\ep^{2}\lambda^{2}(t)
	\quad\quad
	\forall t\geq 0.
	\label{est:z-main-1}
\end{equation}

On the other hand, the estimate from below in (\ref{th:par-Phi}) 
implies that
\begin{equation}
	e^{2\eta\Phi(t)}\geq k_{3}(1+t)^{\eta}
	\quad\quad
	\forall t\geq 0.
	\label{est:z-main-2}
\end{equation}

From (\ref{est:z-main-1}) and (\ref{est:z-main-2}) we easily conclude 
(\ref{th:roep}).

Thus we can limit ourselves to show that $z_{\ep}(t)$ satisfies the
assumptions of Lemma~\ref{lemma:main-ODE}. To this end, we first 
observe that $\roep(t)$ is the solution of the first order equation
$$\roep'(t)=-c(t)A\roep(t)-(\cep(t)-c(t))A\uep(t)-
\ep\uep''(t),$$
with initial condition $\roep(0)=0$.
Therefore we have that $z_{\ep}(0)=0$, and
\begin{eqnarray}
	z_{\ep}'(t) & = & \delta\nu^{2}|A^{1/2}\ul(t)|^{2\gamma}
	\cdot|\roep(t)|^{2}\cdot
	e^{2\eta\Phi(t)}-c(t)\cdot|A^{1/2}\roep(t)|^{2}\cdot e^{2\eta\Phi(t)}
	\nonumber  \\
	\noalign{\vspace{0.5ex}}
	 &  & -[\cep(t)-c(t)]\cdot\langle\roep(t),A\uep(t)\rangle\cdot
	 e^{2\eta\Phi(t)}-\ep\langle\uep''(t),\roep(t)\rangle\cdot
	 e^{2\eta\Phi(t)}
	\nonumber  \\
	 & =: & L_{1}+L_{2}+L_{3}+L_{4}.
	\label{est:z-L1-L4}
\end{eqnarray}

The term $L_{1}$ is $z_{\ep}(t)$ times a coefficient which behaves
like $(1+t)^{-1}$, hence whose integral is divergent.  Therefore this
term alone would prevent~(\ref{est:z-main-1}) from being true.  Thus
the idea is to cancel out $L_{1}$ by means of $L_{2}$ and one of the
terms arising from the expansion of $L_{3}$.  In the following
paragraphs we carry out this program.

\subparagraph{\textmd{\emph{Estimate of $L_{1}$ and $L_{2}$}}}

We claim that
\begin{equation}
	L_{1}+L_{2}\leq 2\gamma|A^{1/2}\ul(t)|^{2\gamma}\cdot
	|A^{1/2}\roepl(t)|^{2}\cdot e^{2\eta\Phi(t)}.
	\label{est:L1+L2}
\end{equation}

Since $|A^{1/2}\roepl(t)|^{2}=\nu^{2}|\roepl(t)|^{2}$, from
(\ref{defn:delta}) we have that
\begin{eqnarray}
	L_{1} & = & \delta|A^{1/2}\ul|^{2\gamma}\cdot
	\nu^{2}|\roepl|^{2}\cdot e^{2\eta\Phi}+
	\delta\nu^{2}|A^{1/2}\ul|^{2\gamma}\cdot
	|\roeph|^{2}\cdot e^{2\eta\Phi}
	\nonumber  \\
	\noalign{\vspace{0.5ex}}
	 & \leq & (2\gamma+1)|A^{1/2}\ul|^{2\gamma}\cdot
	|A^{1/2}\roepl|^{2}\cdot e^{2\eta\Phi}+
	\delta_{0}\nu^{2}|A^{1/2}\ul|^{2\gamma}\cdot
	|\roeph|^{2}\cdot e^{2\eta\Phi}.
	\label{est:L1}
\end{eqnarray}

In order to estimate $L_{2}$, we observe that
$c(t)=|A^{1/2}u(t)|^{2\gamma}\geq|A^{1/2}\ul(t)|^{2\gamma}$, and we
exploit (\ref{coercivity}) with $\mu=\delta_{0}\nu^{2}$ to deduce that
$$|A^{1/2}\roep(t)|^{2}=
|A^{1/2}\roepl(t)|^{2}+|A^{1/2}\roeph(t)|^{2}\geq
|A^{1/2}\roepl(t)|^{2}+\delta_{0}\nu^{2}|\roeph(t)|^{2},$$

It follows that
\begin{equation}
	L_{2}\leq-|A^{1/2}\ul|^{2\gamma}\cdot
	|A^{1/2}\roepl|^{2}\cdot e^{2\eta\Phi}-
	\delta_{0}\nu^{2}|A^{1/2}\ul|^{2\gamma}\cdot|\roeph|^{2}\cdot
	e^{2\eta\Phi}.
	\label{est:L2}
\end{equation}

Adding (\ref{est:L1}) and (\ref{est:L2}) we obtain (\ref{est:L1+L2}).

\subparagraph{\textmd{\emph{Estimate of $L_{3}$}}}

We claim that
\begin{eqnarray}
	L_{3} & \leq & -2\gamma|A^{1/2}\ul(t)|^{2\gamma}\cdot
	|A^{1/2}\roepl(t)|^{2}\cdot e^{2\eta\Phi(t)}
	\nonumber  \\
	 &  & +k_{4}\frac{z_{\ep}(t)}{(1+t)^{1+(\delta-1)/(2\gamma)}}
	 +k_{5}\frac{[z_{\ep}(t)]^{3/2}}{(1+t)^{1+(\delta-1)/(2\gamma)}}.
	\label{est:L3}
\end{eqnarray}

We point out that the first term cancels out the right-hand side of 
(\ref{est:L1+L2}).

In order to prove (\ref{est:L3}), we set
$$R(t):=|A^{1/2}\uep(t)|^{2\gamma}-|A^{1/2}u(t)|^{2\gamma}-
\gamma|A^{1/2}u(t)|^{2(\gamma-1)}\left(
|A^{1/2}\uep(t)|^{2}-|A^{1/2}u(t)|^{2}\right),$$
so that
\begin{equation}
	\cep(t)-c(t)=R(t)+\gamma|A^{1/2}u(t)|^{2(\gamma-1)}\left(
	|A^{1/2}\uep(t)|^{2}-|A^{1/2}u(t)|^{2}\right).
	\label{est:R(t)-1}
\end{equation}

Now we observe that
\begin{eqnarray}
	|A^{1/2}\uep(t)|^{2}-|A^{1/2}u(t)|^{2} & = & 
	\langle\uep(t)-u(t),A\uep(t)+Au(t)\rangle
	\nonumber  \\
	 & = & \langle\roep(t),A\roep(t)+2Au(t)\rangle
	\nonumber  \\
	 & = & |A^{1/2}\roep(t)|^{2}+2\langle\roep(t),Au(t)\rangle,
	\label{est:R(t)-2}
\end{eqnarray}
and
\begin{equation}
	\langle\roep(t),A\uep(t)\rangle=
	\langle\roep(t),A\roep(t)+Au(t)\rangle=
	|A^{1/2}\roep(t)|^{2}+\langle\roep(t),Au(t)\rangle.
	\label{est:R(t)-3}
\end{equation}

Thus from (\ref{est:R(t)-1}), (\ref{est:R(t)-2}), and
(\ref{est:R(t)-3}) we deduce that
\begin{eqnarray*}
	L_{3} & = & -R(t)\cdot\langle\roep(t),A\uep(t)\rangle\cdot
	e^{2\eta\Phi(t)}- 
	\gamma|A^{1/2}u(t)|^{2(\gamma-1)}\cdot e^{2\eta\Phi(t)}\times  \\
	\noalign{\vspace{0.5ex}}
	 &  & \times
	 \left(|A^{1/2}\roep(t)|^{2}+2\langle\roep(t),Au(t)\rangle\right)\cdot
	 \left(|A^{1/2}\roep(t)|^{2}+\langle\roep(t),Au(t)\rangle\right).
\end{eqnarray*}

Neglecting the negative term with $|A^{1/2}\roep(t)|^{4}$, we obtain 
that
\begin{eqnarray}
	L_{3} & \leq & |R(t)|\cdot|\roep(t)|\cdot|A\uep(t)|\cdot e^{2\eta\Phi(t)}
	\nonumber  \\
	\noalign{\vspace{0.5ex}}
	 &  & -3\gamma|A^{1/2}u(t)|^{2(\gamma-1)}\cdot|A^{1/2}\roep(t)|^{2}
	 \cdot\langle\roep(t),Au(t)\rangle\cdot e^{2\eta\Phi(t)}
	\nonumber  \\
	\noalign{\vspace{0.5ex}}
	 &  & -2\gamma|A^{1/2}u(t)|^{2(\gamma-1)}\cdot
	 \langle\roep(t),Au(t)\rangle^{2}\cdot e^{2\eta\Phi(t)}
	\nonumber  \\
	 & =: & L_{3,1}+L_{3,2}+L_{3,3}.
	\nonumber
\end{eqnarray}

Now we claim that
\begin{eqnarray}
	& \displaystyle{L_{3,1} \leq
	k_{6}\frac{[z_{\ep}(t)]^{3/2}}{(1+t)^{1+(\delta-1)/(2\gamma)}},} &
	\label{est:L31} \\
	\noalign{\vspace{0.5ex}}
	& \displaystyle{L_{3,2} \leq k_{7}
	\frac{[z_{\ep}(t)]^{3/2}}{(1+t)^{1+(\delta-1)/(2\gamma)}}+ k_{8}
	\frac{z_{\ep}(t)}{(1+t)^{1+(\delta-1)/(2\gamma)}},} &
	\label{est:L32}  \\
	\noalign{\vspace{0.5ex}}
	& \displaystyle{L_{3,3} \leq -2\gamma|A^{1/2}\ul(t)|^{2\gamma}\cdot
	|A^{1/2}\roepl(t)|^{2}\cdot e^{2\eta\Phi(t)}+
	k_{9}\frac{z_{\ep}(t)}{(1+t)^{1+(\delta-1)/(2\gamma)}},} &
	\label{est:L33}
\end{eqnarray}
from which (\ref{est:L3}) follows directly. The proof of 
(\ref{est:L31}) through (\ref{est:L33}) is the content of the next 
three paragraphs.

\subparagraph{\textmd{\emph{Estimate of $L_{3,1}$}}}

From the second order Taylor's expansion of the function
$\sigma^{\gamma}$ it follows that
$$\left|y^{\gamma}-x^{\gamma}-\gamma x^{\gamma-1}(y-x)\right|\leq
\frac{\gamma(\gamma-1)}{2}\max\left\{x^{\gamma-2},y^{\gamma-2}\right\}
(y-x)^{2} \quad\quad
\forall x\geq 0,\ \forall y\geq 0.$$

Setting $x:=|A^{1/2}u(t)|^{2}$ and $y:=|A^{1/2}\uep(t)|^{2}$, we
obtain that
\begin{equation}
	|R(t)|\leq k_{10}\max\left\{
	|A^{1/2}u|^{2(\gamma-2)},|A^{1/2}\uep|^{2(\gamma-2)}\right\}
	\cdot\left(|A^{1/2}\uep|^{2}-|A^{1/2}u|^{2}\right)^{2}.
	\label{est:R-1}
\end{equation}

Now from (\ref{th:dec-par}) with $j=1$ and (\ref{est:auq}) we have 
that
\begin{equation}
	\max\left\{
	|A^{1/2}u(t)|^{2(\gamma-2)},|A^{1/2}\uep(t)|^{2(\gamma-2)}
	\right\}\leq\frac{k_{11}}{(1+t)^{1-2/\gamma}}
	\label{est:R-2}
\end{equation}
(note that in (\ref{th:dec-par}) and (\ref{est:auq}) we need both the 
estimates from below and the estimates from above because we ignore 
the sign of $\gamma-2$). From (\ref{th:dec-par}) with $j=2$ and 
(\ref{est:au}) we have that
\begin{eqnarray}
	\left(|A^{1/2}\uep(t)|^{2}-|A^{1/2}u(t)|^{2}\right)^{2} & = & 
	\langle\roep(t),Au(t)+A\uep(t)\rangle^{2}
	\nonumber  \\
	 & \leq &|\roep(t)|^{2}\cdot
	 2\left(|Au(t)|^{2}+|A\uep(t)|^{2}\right)
	\nonumber  \\
	 & \leq & k_{12}\frac{|\roep(t)|^{2}}{(1+t)^{1/\gamma}}.
	\label{est:R-3}
\end{eqnarray}

From (\ref{est:R-1}), (\ref{est:R-2}), and (\ref{est:R-3}) it follows that
$$|R(t)|\leq k_{13}\frac{|\roep(t)|^{2}}{(1+t)^{1-1/\gamma}},$$
hence
$$L_{3,1}\leq k_{13}\frac{|\roep(t)|^{3}}{(1+t)^{1-1/\gamma}}\cdot
|A\uep(t)|\cdot e^{2\eta\Phi(t)}\leq
k_{14}\frac{|z(t)|^{3/2}}{(1+t)^{1-1/\gamma}}\cdot
|A\uep(t)|\cdot e^{-\eta\Phi(t)}.$$

The last two terms can be easily estimated exploiting (\ref{est:au})
and the estimate from below in (\ref{th:par-Phi}).  We thus
obtain~(\ref{est:L31}).

\subparagraph{\textmd{\emph{Estimate of $L_{3,2}$}}}

Let us begin by remarking that
\begin{equation}
	L_{3,2}\leq k_{15}|A^{1/2}u(t)|^{2(\gamma-1)}\cdot|A^{1/2}\roep(t)|^{2}
	\cdot|\roep(t)|\cdot|Au(t)|\cdot e^{2\eta\Phi(t)}.
	\label{est:L32-1}
\end{equation}

The first and fourth term can be estimated exploiting 
(\ref{th:dec-par}) with $j=1$ and $j=2$. We obtain that
\begin{equation}
	|A^{1/2}u(t)|^{2(\gamma-1)}\cdot|Au(t)|\leq
	\frac{k_{16}}{(1+t)^{1-1/(2\gamma)}}.
	\label{est:L32-2}
\end{equation}

For the second term we have that
\begin{eqnarray*}
	|A^{1/2}\roep(t)|^{2} & = &
	|A^{1/2}\roepl(t)|^{2}+|A^{1/2}\roeph(t)|^{2} \\
	\noalign{\vspace{0.5ex}}
	 & \leq & \nu^{2}|\roepl(t)|^{2}+|\roeph(t)|\cdot|A\roeph(t)|  \\
	\noalign{\vspace{0.5ex}}
	 & \leq & \nu^{2}|\roep(t)|^{2}+|\roep(t)|\cdot\left(
	 |A\ueph(t)|+|A\uh(t)|\right).
\end{eqnarray*}

The last two terms can be controlled using our estimates for high
frequencies.  From (\ref{th:mist-Auh}) with
$\mu=\delta_{0}\nu^{2}$ and (\ref{th:par-gain}) we obtain that
\begin{equation}
	|A^{1/2}\roep(t)|^{2}\leq \nu^{2}|\roep(t)|^{2}+
	k_{17}\frac{|\roep(t)|}{(1+t)^{\delta_{0}/(2\gamma)}}.
	\label{est:L32-3}
\end{equation}

From (\ref{est:L32-1}) through (\ref{est:L32-3}) it follows that
\begin{eqnarray*}
	L_{3,2} & \leq & k_{18}
	\frac{|\roep(t)|^{3}\cdot e^{2\eta\Phi(t)}}{(1+t)^{1-1/(2\gamma)}}+
	k_{19}
	\frac{|\roep(t)|^{2}\cdot 
	e^{2\eta\Phi(t)}}{(1+t)^{1+(\delta_{0}-1)/(2\gamma)}}  \\
	\noalign{\vspace{0.5ex}}
	 & \leq & k_{20}
	\frac{[z_{\ep}(t)]^{3/2}\cdot e^{-\eta\Phi(t)}}{(1+t)^{1-1/(2\gamma)}}+
	k_{21}
	\frac{z_{\ep}(t)}{(1+t)^{1+(\delta-1)/(2\gamma)}}.  
\end{eqnarray*}

Exploiting the estimate from below in (\ref{th:par-Phi}) we easily
obtain~(\ref{est:L32}).

\subparagraph{\textmd{\emph{Estimate of $L_{3,3}$}}}

First of all we have that
\begin{eqnarray*}
	\langle\roep(t),Au(t)\rangle^{2} & = & \left(\strut
	\langle\roepl(t),A\ul(t)\rangle+
	\langle\roeph(t),A\uh(t)\rangle\right)^{2}  \\
	 & \geq & \langle\roepl(t),A\ul(t)\rangle^{2}+
	 2\langle\roepl(t),A\ul(t)\rangle\cdot
	 \langle\roeph(t),A\uh(t)\rangle \\
	 & \geq & \langle A^{1/2}\roepl(t),A^{1/2}\ul(t)\rangle^{2}-
	 2\,|\roep(t)|^{2}\cdot|A\ul(t)|\cdot|A\uh(t)|,
\end{eqnarray*}
hence
\begin{eqnarray}
	L_{3,3} & \leq & -2\gamma|A^{1/2}u(t)|^{2(\gamma-1)}\cdot\langle
	A^{1/2}\roepl(t),A^{1/2}\ul(t)\rangle^{2}\cdot e^{2\eta\Phi(t)}
	\nonumber \\
	\noalign{\vspace{0.5ex}}
	 &  & +4\gamma|A^{1/2}u(t)|^{2(\gamma-1)}\cdot
	 |\roep(t)|^{2}\cdot|A\ul(t)|\cdot|A\uh(t)|\cdot e^{2\eta\Phi(t)}
	\nonumber  \\
	 & =: & L_{3,3,1}+L_{3,3,2}.
	 \label{est:L33-split}
\end{eqnarray}

Since $u_{0,\ell}$ and $u_{1,\ell}$ are multiples of an eigenvector of
$A$, it is easy to see that both $\uepl(t)$ and $\ul(t)$ are multiples
of the same eigenvector, and the same for $\roepl(t)$.  Therefore the
vectors $A^{1/2}\roepl(t)$ and $A^{1/2}\ul(t)$ are parallel, hence the
square of their scalar product is equal to the square of the product
of their norms (this is the point where the last condition in the
$(\nu,\delta_{0})$-assumption plays a crucial role).  It follows that
\begin{eqnarray}
	L_{3,3,1} & = & -2\gamma|A^{1/2}u(t)|^{2(\gamma-1)}\cdot
	|A^{1/2}\roepl(t)|^{2}\cdot|A^{1/2}\ul(t)|^{2}\cdot e^{2\eta\Phi(t)}
	\nonumber  \\
	 & \leq & -2\gamma|A^{1/2}\ul(t)|^{2\gamma}\cdot
	|A^{1/2}\roepl(t)|^{2}\cdot e^{2\eta\Phi(t)}.
	\label{est:L331}
\end{eqnarray}

On the other hand, exploiting (\ref{th:dec-par}) with $j=1$ and 
$j=2$, and (\ref{th:par-gain}), we have that
\begin{eqnarray}
	L_{3,3,2} & \leq & k_{22}|A^{1/2}u(t)|^{2(\gamma-1)}\cdot
	|Au(t)|\cdot|A\uh(t)|\cdot|\roep(t)|^{2}\cdot e^{2\eta\Phi(t)}
	\nonumber  \\
	\noalign{\vspace{0.5ex}}
	 & \leq & k_{23}\frac{1}{(1+t)^{1-1/\gamma}}\cdot
	 \frac{1}{(1+t)^{1/(2\gamma)}}\cdot
	 \frac{1}{(1+t)^{\delta_{0}/(2\gamma)}}\cdot
	 z_{\ep}(t).
	\label{est:L332}
\end{eqnarray}

Plugging (\ref{est:L331}) and (\ref{est:L332}) into
(\ref{est:L33-split}), and recalling that $\delta\leq\delta_{0}$, we
obtain (\ref{est:L33}).

\subparagraph{\textmd{\emph{Estimate of $L_{4}$}}}

Let us fix $\mu:=8\gamma\nu^{2}$.  Splitting components corresponding
to low and high frequencies with respect to $\mu$, we have that
\begin{equation}
	L_{4}=-\ep\langle\ueplm''(t),\roeplm(t)\rangle\cdot 
	e^{2\eta\Phi(t)}-
	\ep\langle\uephh''(t),\roephh(t)\rangle\cdot e^{2\eta\Phi(t)}.
	\label{est:L4-sum}
\end{equation}

Now we claim that
\begin{eqnarray}
	\left|\ep\langle\ueplm''(t),\roeplm(t)\rangle\cdot
	e^{2\eta\Phi(t)}\right| & \leq & \left(
	\frac{\ep}{(e+t)^{2-(\delta-1)/(2\gamma)}}+
	(1+t)^{\eta/2}e^{-t/\ep}\right)\times \nonumber \\
	 &  & \times\left(\frac{z_{\ep}(t)}{\ep\lambda(t)}+
	k_{24}\ep\lambda(t)\right).
	\label{est:L41}
\end{eqnarray}

Indeed from Proposition~\ref{prop:mist-u''} we have that
\begin{eqnarray*}
	\left|\langle\ueplm''(t),\roeplm(t)\rangle\right| 
	 & \leq & |\ueplm''(t)|\cdot|\roeplm(t)|  \\
	 & \leq & \left(|\ueplm''(t)-\teteplm''(t)|+
	 |\teteplm''(t)|\right)\cdot|\roeplm(t)|  \\
	 & \leq & k_{25}\left(
	\frac{1}{(1+t)^{2+1/(2\gamma)}}+
	\frac{1}{\ep}e^{-t/\ep}\right)\cdot|\roep(t)|,
\end{eqnarray*}
so that the estimate from above in (\ref{th:par-Phi}) implies that
\begin{eqnarray*}
	\ep\left|\langle\ueplm''(t),\roeplm(t)\rangle\right| 
	e^{2\eta\Phi(t)} & \leq & k_{26}\ep\left(
	\frac{1}{(1+t)^{2+1/(2\gamma)}}+
	\frac{1}{\ep}e^{-t/\ep}\right)\sqrt{z_{\ep}(t)}\cdot e^{\eta\Phi(t)}  \\
	 & \leq & k_{27}\left(
	\frac{\ep}{(e+t)^{2+1/(2\gamma)}}+
	e^{-t/\ep}\right)\sqrt{z_{\ep}(t)}\cdot(1+t)^{\eta/2}  \\
	 & = & k_{27}\left(
	\frac{\ep}{(e+t)^{2-(\delta-1)/(2\gamma)}}+
	(1+t)^{\eta/2}e^{-t/\ep}\right)\sqrt{z_{\ep}(t)}.
\end{eqnarray*}

Since
$$\sqrt{z_{\ep}(t)}\leq\frac{z_{\ep}(t)}{k_{27}\ep\lambda(t)}+
k_{27}\ep\lambda(t),$$
we have proved (\ref{est:L41}).

\subparagraph{\textmd{\emph{Checking the assumptions of 
Lemma~\ref{lemma:main-ODE}}}}

Plugging (\ref{est:L1+L2}), (\ref{est:L3}), (\ref{est:L4-sum}), and
(\ref{est:L41}) into (\ref{est:z-L1-L4}), we obtain that $z_{\ep}(t)$ 
satisfies a differential inequality such as (\ref{eqn:zep}) with
\begin{eqnarray*}
	 & \displaystyle{\psi_{1}(\ep,t):=
	 \frac{1}{(e+t)^{2-(\delta-1)/(2\gamma)}}\cdot
	 \frac{1}{\lambda(t)}+
	 \frac{1}{\ep}(1+t)^{\delta/(2\gamma)}
	 \frac{1}{\lambda(t)}e^{-t/\ep}+
	 \frac{k_{4}}{(1+t)^{1+(\delta-1)/(2\gamma)}},} &   \\
	 & \displaystyle{\psi_{2}(\ep,t):=
	 \frac{k_{5}}{(1+t)^{1+(\delta-1)/(2\gamma)}},} &   \\
	 & \displaystyle{\psi_{3}(\ep,t):=
	 k_{24}\ep^{2}\frac{\lambda(t)}{(e+t)^{2-(\delta-1)/(2\gamma)}}
	 +k_{24}\ep\lambda(t)(1+t)^{\delta/(2\gamma)} e^{-t/\ep}, } &  \\
	 & \psi_{4}(\ep,t):=
	 -\ep\langle\uephh''(t),\roephh(t)\rangle\cdot e^{2\eta\Phi(t)}. & 
\end{eqnarray*}

In order to apply Lemma~\ref{lemma:main-ODE} we have to check 
assumptions (\ref{hp:psi0}) through (\ref{hp:psi4}).

Assumption (\ref{hp:psi0}) is trivial.

Let us prove the estimate from below in (\ref{hp:psi1}). If 
$\delta<2\gamma+1$ this is trivial because $\lambda(t)\equiv 1$. If 
$\delta=2\gamma+1$, then we have that
\begin{equation}
	2-\frac{\delta-1}{2\gamma}=1.
	\label{esponente}
\end{equation}

Limiting ourselves to the first term in the expression of 
$\psi_{1}(\ep,t)$, we have therefore that
$$\int_{0}^{t}\psi_{1}(\ep,s)\,ds\geq
\int_{0}^{t}\frac{ds}{(e+s)\lambda(s)}=
\int_{0}^{t}\frac{ds}{(e+s)\log(e+s)}=
\log(\log(e+t))=\log\lambda(t).$$

Let us prove the estimate from above in (\ref{hp:psi1}). Since 
$\delta>1$ we have that
$$\int_{0}^{+\infty}\frac{1}{(1+s)^{1+(\delta-1)/(2\gamma)}}\,ds\leq 
k_{28},$$
and this settles the integral of the third term in the definition of
$\psi_{1}(\ep,t)$. For the integral of the second term, we exploit 
that $\lambda(t)\geq 1$, and with the variable change $\sigma=s/\ep$ 
we obtain that
\begin{equation}
	\frac{1}{\ep}\int_{0}^{+\infty}(1+s)^{\delta/(2\gamma)}
	\frac{1}{\lambda(s)}e^{-s/\ep}\,ds\leq
	\int_{0}^{+\infty}(1+\ep_{0}\sigma)^{\delta/(2\gamma)}
	e^{-\sigma}\,d\sigma\leq k_{29}.
	\label{est:int-exp}
\end{equation}

It remains to estimate the integral of the first term.  If
$\delta<2\gamma+1$ we have that
$$\int_{0}^{+\infty}\frac{1}{(e+s)^{2-(\delta-1)/(2\gamma)}}
\cdot\frac{1}{\lambda(s)}\,ds=
\int_{0}^{+\infty}\frac{1}{(e+s)^{2-(\delta-1)/(2\gamma)}}\,ds
\leq k_{30}.$$

If $\delta=2\gamma+1$, then by (\ref{esponente}) we have that
$$\int_{0}^{t}\frac{1}{(e+s)^{2-(\delta-1)/(2\gamma)}}
\cdot\frac{1}{\lambda(s)}\,ds=
\int_{0}^{t}\frac{ds}{(e+s)\log(e+s)}=
\log(\log(e+t))=\log\lambda(t).$$

In both cases we have proved (\ref{hp:psi1}).

Let us consider now (\ref{hp:psi2}). Since $\delta>1$ we have that
$$\int_{0}^{+\infty}|\psi_{2}(\ep,s)|\cdot\lambda^{3}(s)\,ds\leq
k_{5}\int_{0}^{+\infty}
\frac{\log^{3}(e+s)}{(1+s)^{1+(\delta-1)/(2\gamma)}}\,ds\leq k_{31},$$
which proves (\ref{hp:psi2}).

In order to prove (\ref{hp:psi3}), we consider the integral
$$\int_{0}^{t}\frac{|\psi_{3}(\ep,s)|}{\lambda(s)}\,ds=
k_{24}\ep^{2}\int_{0}^{t}\frac{ds}{(e+s)^{2-(\delta-1)/(2\gamma)}}+
k_{24}\ep\int_{0}^{t}(1+s)^{\delta/(2\gamma)}
e^{-s/\ep}\,ds.$$

The second integral can be estimated as in (\ref{est:int-exp}). The 
first integral is less than a constant if $\delta<2\gamma+1$, and equal to 
$\log(e+t)=\lambda(t)$ if $\delta=2\gamma+1$. In both cases this 
proves (\ref{hp:psi3}).

It remains to prove (\ref{hp:psi4}), and this is the content of the 
last paragraph.

\subparagraph{\textmd{\emph{Estimate of the integral of 
$\psi_{4}(\ep,t)$}}}

We have to prove that
\begin{equation}
	\left|\int_{0}^{t}
	\ep\langle\uephh''(s),\roephh(s)\rangle\cdot
	e^{2\eta\Phi(s)}\,ds
	\right|\leq k_{32}\ep^{2}
	\quad\quad
	\forall t\geq 0.
	\label{est:L42}
\end{equation}

To this end, we first integrate by parts and we obtain that
\begin{eqnarray}
	\lefteqn{\hspace{-5em}
	-\int_{0}^{t} \ep\langle\uephh''(s),\roephh(s)\rangle\cdot
	e^{2\eta\Phi(s)}\,ds\  =\  -\ep\langle\uephh'(t),\roephh(t)\rangle
	\cdot e^{2\eta\Phi(t)}}
	\nonumber  \\
	\quad\quad\quad & & +\int_{0}^{t}
	\ep\langle\uephh'(s),\roephh'(s)\rangle\cdot e^{2\eta\Phi(s)}\,ds
	\nonumber \\
	 &  & +2\delta\nu^{2}\int_{0}^{t} \ep\langle\uephh'(s),\roephh(s)\rangle
	 \cdot |A^{1/2}\ul(s)|^{2\gamma}\cdot e^{2\eta\Phi(s)}\,ds
	\nonumber  \\
	 & =: & I_{1}+I_{2}+I_{3}.
	\label{est:L42-I}
\end{eqnarray}

Now from (\ref{th:mist-u'}) with $\mu=8\gamma\nu^{2}$ we have that
\begin{equation}
	|\uephh'(t)|\leq\frac{k_{33}}{(1+t)^{5}}
	\quad\quad
	\forall t\geq 0,
	\label{est:uephh'}
\end{equation}
and from Theorem~\ref{thm:hyst-sp} we have that
\begin{equation}
	|\roep(t)|^{2}\leq k_{34}\ep^{2}
	\quad\quad
	\forall t\geq 0,
	\label{est:roep-C}
\end{equation}
\begin{equation}
	\int_{0}^{+\infty}(1+s)|\rep'(s)|^{2}\,ds\leq k_{35}\ep^{2}.
	\label{est:rep-C}
\end{equation}

Exploiting (\ref{est:uephh'}), (\ref{est:roep-C}), the estimate from
above in (\ref{th:par-Phi}), and the fact that $\eta\leq
2+1/\gamma\leq 5$, we have that
\begin{equation}
	|I_{1}|\leq\ep|\uephh'(t)|\cdot|\roep(t)|\cdot e^{2\eta\Phi(t)}
	\leq k_{36}\frac{1}{(1+t)^{5}}\cdot\ep^{2}\cdot(1+t)^{\eta}
	\leq k_{37}\ep^{2}.
	\label{est:L42-I1}
\end{equation}

In order to estimate $|I_{2}|$, we first observe that
\begin{equation}
	|\roephh'(t)|\leq|\roep'(t)|\leq|\rep'(t)|+|\tetep'(t)|\leq
	|\rep'(t)|+k_{38}e^{-t/\ep}.
	\label{est:roephh'}
\end{equation}

From (\ref{est:uephh'}), the estimate from above in
(\ref{th:par-Phi}), and the fact that $\eta\leq 2+1/\gamma\leq 3$, we
have that
\begin{equation}
	|\uephh'(t)|\cdot e^{2\eta\Phi(t)}\leq
	k_{39}\frac{1}{(1+t)^{5}}\cdot (1+t)^{\eta}
	\leq k_{39}\frac{1}{(1+t)^{2}}.
	\label{est:uephh'-bis}
\end{equation}

Thanks to (\ref{est:roephh'}) and (\ref{est:uephh'-bis}) we obtain that
\begin{eqnarray*}
	\left|\ep\langle\uephh'(t),\roephh'(t)\rangle\cdot
	e^{2\eta\Phi(t)}\right| & \leq & \ep|\uephh'(t)|\cdot|\roephh'(t)|
	\cdot e^{2\eta\Phi(t)} \\
	 & \leq & k_{39}\ep\frac{|\rep'(t)|}{(1+t)^{2}}+
	 k_{40}\ep\frac{e^{-t/\ep}}{(1+t)^{2}}\\
	 & \leq & (1+t)|\rep'(t)|^{2}+k_{41}\ep^{2}\frac{1}{(1+t)^{5}}
	 +k_{40}\ep e^{-t/\ep}.
\end{eqnarray*}

Integrating in $[0,t]$, and exploiting (\ref{est:rep-C}), we deduce
that
\begin{equation}
	|I_{2}|\leq
	\int_{0}^{t}(1+s)|\rep'(s)|^{2}\,ds+
	k_{41}\ep^{2}\int_{0}^{t}\frac{1}{(1+s)^{5}}\,ds
	 +k_{40}\ep \int_{0}^{t}e^{-s/\ep}\,ds\leq k_{42}\ep^{2}.
	\label{est:L42-I2}
\end{equation}

Finally, from (\ref{est:uephh'}), (\ref{est:roep-C}),
(\ref{th:dec-par}) with $j=1$, and the estimate from above in
(\ref{th:par-Phi}), we obtain that
\begin{eqnarray*}
	|\langle\uephh'(t),\roephh(t)\rangle|\cdot
	|A^{1/2}u_{l}(t)|^{2\gamma}\cdot e^{2\eta\Phi(s)} & \leq &
	|\uephh'(t)|\cdot|\roep(t)|\cdot
	|A^{1/2}u(t)|^{2\gamma}\cdot e^{2\eta\Phi(t)}\\
	 & \leq & k_{43}\frac{1}{(1+t)^{5}}\cdot\ep\cdot
	 \frac{1}{1+t}\cdot(1+t)^{\eta}.
\end{eqnarray*}

Since $\eta\leq 3$, we conclude that
\begin{equation}
	|I_{3}|\leq k_{44}\ep^{2}
	\int_{0}^{t}\frac{1}{(1+s)^{6-\eta}}\,ds\leq k_{45}\ep^{2}.
	\label{est:L42-I3}
\end{equation}

Plugging (\ref{est:L42-I1}), (\ref{est:L42-I2}), and
(\ref{est:L42-I3}) into (\ref{est:L42-I}) we obtain (\ref{est:L42}).

This completes the proof of~(\ref{th:roep}).\qed

\subsection{Estimates on linear equations}

Let us define $\cep(t)$ and $c(t)$ as in (\ref{defn:c-cep}), and let us
set
\begin{equation}
	\gep(t):=-(\cep(t)-c(t))Au(t)-\ep u''(t).
	\label{defn:gep}
\end{equation}

Then it is easy to see that $\roep(t)$ is the solution of the linear
equation
\begin{equation}
	\ep\roep''(t)+\roep'(t)+\cep(t)A\roep(t)=\gep(t),
	\label{eqn:roep-lin}
\end{equation}
with initial data
\begin{equation}
	\roep(0)=0,
	\quad\quad
	\roep'(0)=w_{0},
	\label{eqn:roep-data}
\end{equation}
while $\rep(t)$ is the solution of the linear equation
\begin{equation}
	\ep\rep''(t)+\rep'(t)+\cep(t)A\roep(t)=\gep(t),
	\label{eqn:rep-lin}
\end{equation}
with initial data
\begin{equation}
	\rep(0)=0,
	\quad\quad
	\rep'(0)=0.
	\label{eqn:rep-data}
\end{equation}

In this section we forget that $\gep(t)$ is given by (\ref{defn:gep}), 
and that $\cep(t)$ and $c(t)$ are given by (\ref{defn:c-cep}). We just regard 
$\roep(t)$ and $\rep(t)$ as solutions of the corresponding linear 
equations (which implies also that $\roep(t)=\rep(t)+\tetep(t)$, 
where $\tetep(t)$ given by (\ref{eqn:tetep})).

We assume that the coefficient $\cep:[0,+\infty)\to(0,+\infty)$ is of
class $C^{1}$ and such that
\begin{equation}
	\frac{M_{1}}{1+t}\leq\cep(t)\leq\frac{M_{2}}{1+t}
	\quad\quad
	\forall t\geq 0,
	\label{hp:cep-1}
\end{equation}
\begin{equation}
	\frac{|\cep'(t)|}{\cep(t)}\leq\frac{M_{3}}{1+t}
	\quad\quad
	\forall t\geq 0.
	\label{hp:cep-2}
\end{equation}

We assume that the forcing term 
$\gep:[0,+\infty)\to H$ is continuous and such that
\begin{eqnarray}
	 & \displaystyle{\int_{0}^{t}(1+s)^{1+2\delta/\gamma}|\gep(s)|^{2}\,ds\leq
	 M_{4}\ep^{2}(1+t)^{\delta/\gamma}\lambda^{2}(t)}
	 \quad\quad\forall t\geq 0, & 
	\label{hp:gep-1}  \\
	 & \displaystyle{\int_{0}^{t}
	 (1+s)^{1+2\delta/\gamma}|A^{1/2}\gep(s)|^{2}\,ds\leq
	 M_{5}\ep^{2}(1+t)^{\delta/\gamma}\lambda^{2}(t)} 
	 \quad\quad\forall t\geq 0,  & 
	\label{hp:gep-2}  \\
	 & \displaystyle{|\gep(t)|^{2}\leq M_{6}\ep^{2}
	 \frac{\lambda^{2}(t)}{(1+t)^{2+\delta/\gamma}}} 
	 \quad\quad\forall t\geq 0,  & 
	\label{hp:gep-3}
\end{eqnarray}
where $\gamma>0$, $\delta>0$, and $\lambda:[0,+\infty)\to[1,+\infty)$
is a continuous nondecreasing function such that
\begin{equation}
	\int_{0}^{t} (1+s)^{-1+\delta/\gamma}\cdot\lambda^{2}(s)\,ds\leq
	M_{7}(1+t)^{\delta/\gamma}\lambda^{2}(t).
	\label{est:int-2}
\end{equation}

These requirements on $\gamma$, $\delta$, $\lambda(t)$ are weaker than 
those in Theorem~\ref{thm:main}.

Under such assumptions we show that an a priori estimate on $\roep(t)$
of the form~(\ref{th:roep}) yields all other estimates on $\roep$,
$\rep$, and their derivatives contained in statements~(1) and~(2) of
Theorem~\ref{thm:main}.

\begin{prop}\label{prop:linear}
	Let $H$ be a Hilbert space, let $A$ be a nonnegative self-adjoint
	(unbounded) operator on $H$ with dense domain, and let
	$\ep_{0}>0$.  For every $\ep\in (0,\ep_{0})$, let $\roep(t)$ and
	$\rep(t)$ be functions in the space (\ref{reg:space}) satisfying
	(\ref{eqn:roep-lin}) through (\ref{eqn:rep-data}).  Let us assume
	that $\cep(t)$, $\gep(t)$, $\gamma$, $\delta$, $\lambda(t)$ satisfy
	conditions (\ref{hp:cep-1}), (\ref{hp:cep-2}), (\ref{hp:gep-1}) 
	and (\ref{est:int-2}) as above.
	
	Let us assume that $\roep(t)$ satisfies the a priori estimate 
	(\ref{th:roep}).
	
	Then the following conclusions hold true.
	\begin{enumerate}
		\renewcommand{\labelenumi}{(\arabic{enumi})}
		\item  If $w_{0}\in\dau$, then all the estimates in 
		statement~(1) of Theorem~\ref{thm:main} hold true.
	
		\item If in addition $w_{0}\in\da$, and $\gep(t)$ satisfies
		also (\ref{hp:gep-2}) and (\ref{hp:gep-3}), then all the
		estimates in statement~(2) of Theorem~\ref{thm:main} hold
		true.
	\end{enumerate}
\end{prop}

\paragraph{\textmd{\emph{Proof}}}

Let us consider the following weighted versions of classical energies
\begin{eqnarray*}
	\Dep(t) & := & (1+t)^{2\dg}\left(
	\ep\langle\roep(t),\roep'(t)\rangle+
	\frac{1}{2}|\roep(t)|^{2}\right), \\
	\Eep(t) & := & (1+t)^{2\dg}\left(
	\ep\frac{|\rep'(t)|^{2}}{\cep(t)}+|A^{1/2}\roep(t)|^{2}\right).
\end{eqnarray*}

Exploiting (\ref{eqn:roep-lin}) and (\ref{eqn:rep-lin}), with some 
computations we obtain that
\begin{equation}
	\Dep'(t) \ =\  \frac{2\dg\Dep(t)}{1+t}+ (1+t)^{2\dg}\left(
	-\cep(t)|A^{1/2}\roep(t)|^{2}+\ep|\roep'(t)|^{2}+
	\langle\gep(t),\roep(t)\rangle\right),
	\label{deriv-Dep}
\end{equation}
\begin{eqnarray}
	\Eep'(t) & = & -(1+t)^{2\dg}\left(
	2+\ep\frac{\cep'(t)}{\cep(t)}-
	\frac{2\dg\ep}{1+t}\right)
	\frac{|\rep'(t)|^{2}}{\cep(t)}+
	2\dg(1+t)^{2\dg-1}|A^{1/2}\roep(t)|^{2}
	\nonumber\\
	 &  & +2(1+t)^{2\dg}\langle A\roep(t),\tetep'(t)\rangle
	 +2(1+t)^{2\dg}\frac{1}{\cep(t)}
	 \langle\gep(t),\rep'(t)\rangle.
	 \label{deriv-Eep}
\end{eqnarray}

\subparagraph{\textmd{\emph{First energy estimate}}}

We prove that
\begin{eqnarray}
	\int_{0}^{t}(1+s)^{2\dg-1}|A^{1/2}\roep(s)|^{2}\,ds & \leq &
	k_{1}\ep\Eep(t)+k_{2}\ep\int_{0}^{t}(1+s)^{2\dg}
	\frac{|\rep'(s)|^{2}}{\cep(s)}\,ds
	\nonumber  \\
	 & & +k_{3}\ep^{2} (1+t)^{\dg}\lambda^{2}(t).
	\label{est:lin-1}
\end{eqnarray}

To this end, from (\ref{deriv-Dep}) we have that
\begin{eqnarray}
	\Dep'(t) & \leq & \frac{2\dg}{1+t}|\Dep(t)|-
	\cep(t)(1+t)^{2\dg}|A^{1/2}\roep(t)|^{2}
	\nonumber  \\
	\noalign{\vspace{0.5ex}}
	 &  & + \ep(1+t)^{2\dg}|\roep'(t)|^{2}+
	(1+t)^{2\dg}|\gep(t)|\cdot|\roep(t)|
	\nonumber  \\
	\noalign{\vspace{0.5ex}}
	 & =: & L_{1}+L_{2}+L_{3}+L_{4}.
	\label{est:deriv-Dep}
\end{eqnarray}

Let us estimate the four terms.  From (\ref{th:roep}) we have that
\begin{eqnarray*}
	L_{1} & \leq & k_{4}(1+t)^{2\dg-1}\left(
	\ep|\roep'(t)|^{2}+|\roep(t)|^{2}\right)\\
	\noalign{\vspace{0.5ex}}
	 & \leq & k_{4}\ep(1+t)^{2\dg}|\roep'(t)|^{2}+
	k_{5}\ep^{2}(1+t)^{\dg-1}\lambda^{2}(t), 
\end{eqnarray*}
hence
\begin{eqnarray}
	L_{1}+L_{3} & \leq & k_{6}\ep(1+t)^{2\dg}|\roep'(t)|^{2}+
	k_{5}\ep^{2}(1+t)^{\dg-1}\lambda^{2}(t)  
	\nonumber  \\
	\noalign{\vspace{0.5ex}}
	 & \leq  & 2k_{6}\ep(1+t)^{2\dg}\left(
	 |\rep'(t)|^{2}+|\tetep'(t)|^{2}\right)+
	k_{5}\ep^{2}(1+t)^{\dg-1}\lambda^{2}(t).
	\label{est:L1-l}
\end{eqnarray}

From (\ref{hp:cep-1}) we deduce that $\cep(t)$ is bounded.  Therefore,
using also the explicit expression (\ref{eqn:tetep}) for $\tetep(t)$,
from (\ref{est:L1-l}) we obtain that
\begin{equation}
	L_{1}+L_{3}\leq k_{7}\ep(1+t)^{2\dg}
	\frac{|\rep'(t)|^{2}}{\cep(t)}+
	k_{8}\ep(1+t)^{2\dg}e^{-2t/\ep}+
	k_{5}\ep^{2}(1+t)^{\dg-1}\lambda^{2}(t).
	\label{est:L1+L3}
\end{equation}

From (\ref{th:roep}) we obtain also that
\begin{eqnarray}
	L_{4} & \leq & \frac{1}{2}(1+t)^{2\dg-1}|\roep(t)|^{2}+
	\frac{1}{2}(1+t)^{2\dg+1}|\gep(t)|^{2}
	\nonumber \\
	 & \leq &
	k_{9}\ep^{2}(1+t)^{\dg-1}\lambda^{2}(t)+
	\frac{1}{2}(1+t)^{2\dg+1}|\gep(t)|^{2}.
	\label{est:L4}
\end{eqnarray}

Plugging (\ref{est:L1+L3}) and (\ref{est:L4}) into (\ref{est:deriv-Dep}),
and integrating in $[0,t]$, we obtain that
\begin{eqnarray}
	\int_{0}^{t}\cep(s)(1+s)^{2\dg}|A^{1/2}\roep(s)|^{2}\,ds
	 & \leq & -\Dep(t)+k_{7}\ep\int_{0}^{t}(1+s)^{2\dg}
	\frac{|\rep'(s)|^{2}}{\cep(s)}\,ds
	\nonumber  \\
	 &  & +k_{8}\ep\int_{0}^{t}(1+s)^{2\dg}e^{-2s/\ep}\,ds
	\nonumber \\
	&  & + k_{10}\ep^{2}\int_{0}^{t} (1+s)^{\dg-1}\lambda^{2}(s)\,ds
	\nonumber \\
	& & + \frac{1}{2}\int_{0}^{t}(1+s)^{2\dg+1}|\gep(s)|^{2}\,ds.
	\label{est:int-Dep}
\end{eqnarray}

Let us estimate some of the terms in the right-hand side.  Exploiting
the fact that $\cep(t)$ is bounded, the explicit formula
(\ref{eqn:tetep}) for $\tetep(t)$, and the fact that $\lambda(t)\geq
1$, for the first term we obtain that
\begin{eqnarray}
	-\Dep(t) & \leq & \frac{1}{2}\ep^{2}(1+t)^{2\dg}|\roep'(t)|^{2}
	\nonumber  \\
	\noalign{\vspace{0.5ex}}
	 & \leq & \ep^{2}(1+t)^{2\dg}
	 \left(|\rep'(t)|^{2}+|\tetep'(t)|^{2}\right)
	\nonumber  \\
	 & \leq & k_{11}\ep^{2}(1+t)^{2\dg}\frac{|\rep'(t)|^{2}}{\cep(t)}
	 +k_{12}\ep^{2}(1+t)^{2\dg}e^{-2t/\ep}
	\nonumber  \\
	 & \leq & k_{11}\ep\Eep(t)+k_{12}\ep^{2}(1+t)^{\dg}\cdot
	 (1+t)^{\dg}e^{-2t/\ep_{0}}
	\nonumber  \\
	\noalign{\vspace{0.5ex}}
	 & \leq &  k_{11}\ep\Eep(t)+k_{13}\ep^{2}(1+t)^{\dg}\lambda^{2}(t).
	\nonumber
\end{eqnarray}

In the third term of (\ref{est:int-Dep}) we make the variable change
$\sigma=s/\ep$, and we obtain that
\begin{equation}
	\int_{0}^{t}(1+s)^{2\dg}e^{-2s/\ep}\,ds\leq
	\ep\int_{0}^{+\infty}(1+\ep_{0}\sigma)^{2\dg}e^{-2\sigma}\,d\sigma
	\leq k_{14}\ep.
	\label{est:int-1}
\end{equation}

Moreover, we estimate the fourth term of (\ref{est:int-Dep}) by means of
(\ref{est:int-2}), and the fifth by means of~(\ref{hp:gep-1}).
Finally, the estimate from below in (\ref{hp:cep-1}) implies that
$$\int_{0}^{t}(1+s)^{2\dg-1}|A^{1/2}\roep(s)|^{2}\,ds\leq
k_{16}\int_{0}^{t}\cep(s)(1+s)^{2\dg}
|A^{1/2}\roep(s)|^{2}\,ds.$$

Plugging all these estimates into (\ref{est:int-Dep}), we obtain
(\ref{est:lin-1}).

\subparagraph{\textmd{\emph{Second energy estimate}}}

We prove that $\roep(t)$ and $\rep(t)$ satisfy all the conclusions of 
statement~(1) of Theorem~\ref{thm:main}.

We begin by estimating some terms in (\ref{deriv-Eep}).  Thanks to
(\ref{hp:cep-2}) we have that
\begin{equation}
	2+\ep\frac{\cep'(t)}{\cep(t)}- 
	\frac{2\dg\ep}{1+t}\geq\frac{3}{2}
	\quad\quad
	\forall t\geq 0
	\label{est:deriv-Eep-1}
\end{equation}
provided that $\ep$ is small enough. Thanks to the estimate from 
below in (\ref{hp:cep-1}) we have that
\begin{eqnarray}
	2(1+t)^{2\dg}\frac{1}{\cep(t)} \langle\gep(t),\rep'(t)\rangle &
	\leq & \frac{1}{2}(1+t)^{2\dg}\frac{|\rep'(t)|^{2}}{\cep(t)}+
	2(1+t)^{2\dg}\frac{|\gep(t)|^{2}}{\cep(t)}
	\nonumber \\
	 & \leq & \frac{1}{2}(1+t)^{2\dg}\frac{|\rep'(t)|^{2}}{\cep(t)}+
	k_{17}(1+t)^{2\dg+1}|\gep(t)|^{2}.
	\label{est:deriv-Eep-2}
\end{eqnarray}

Moreover we have that
\begin{equation}
	\langle A\roep(t),\tetep'(t)\rangle\leq
	|A^{1/2}\roep(t)|\cdot|A^{1/2}\tetep'(t)|\leq
	|A^{1/2}\roep(t)|\cdot k_{18}e^{-t/\ep}.
	\label{est:deriv-Eep-3}
\end{equation}

Plugging (\ref{est:deriv-Eep-1}) through
(\ref{est:deriv-Eep-3}) into (\ref{deriv-Eep}), we obtain that
\begin{eqnarray*}
	\Eep'(t) & \leq & -(1+t)^{2\dg}\frac{|\rep'(t)|^{2}}{\cep(t)}+
	2\dg(1+t)^{2\dg-1}|A^{1/2}\roep(t)|^{2}\\
	 &  & +2k_{18}(1+t)^{2\dg}|A^{1/2}\roep(t)|\cdot e^{-t/\ep}+
	 k_{17}(1+t)^{2\dg+1}|\gep(t)|^{2}.
\end{eqnarray*}

Integrating in $[0,t]$, and exploiting (\ref{est:lin-1}) and 
(\ref{hp:gep-1}), we obtain that
\begin{eqnarray*}
	\lefteqn{\hspace{-4em}
	\left(1-2\dg \max\{k_{1},k_{2}\}\ep\right)\left(
	\Eep(t)+\int_{0}^{t}(1+s)^{2\dg}\frac{|\rep'(s)|^{2}}{\cep(s)}\,ds
	\right)}  \\
	\hspace{3em} & \leq & k_{19}\ep^{2}(1+t)^{\dg}\lambda^{2}(t)+
	 2k_{18}\int_{0}^{t}(1+s)^{2\dg}|A^{1/2}\roep(s)|
	 \cdot e^{-s/\ep}\,ds.
\end{eqnarray*}

If $\ep$ is small enough this means that
\begin{eqnarray*}
	\lefteqn{\hspace{-4em}
	\Eep(t)+\int_{0}^{t}(1+s)^{2\dg}\frac{|\rep'(s)|^{2}}{\cep(s)}\,ds}  \\
	\hspace{3em} & \leq & k_{20}\ep^{2}(1+t)^{\dg}\lambda^{2}(t)+
	 k_{21}\int_{0}^{t}(1+s)^{2\dg}|A^{1/2}\roep(s)|\cdot e^{-s/\ep}\,ds.
\end{eqnarray*}

Let us fix any $T\geq 0$.  The same argument exploited in
(\ref{est:int-1}) gives that for every $t\in[0,T]$ we have that
\begin{eqnarray*}
	\int_{0}^{t}(1+s)^{2\dg}|A^{1/2}\roep(s)|\cdot
	e^{-s/\ep}\,ds & \leq & 
	\sup_{\tau\in[0,T]}(1+\tau)^{\dg}|A^{1/2}\roep(\tau)|\cdot
	\int_{0}^{t}(1+s)^{\dg}e^{-s/\ep}\,ds \\
	 & \leq & \sup_{\tau\in[0,T]}\sqrt{\Eep(\tau)}\cdot k_{22}\ep  \\
	 & \leq & \frac{1}{2k_{21}}\sup_{\tau\in[0,T]}\Eep(\tau)+k_{23}\ep^{2},
\end{eqnarray*}
hence
\begin{equation}
	\Eep(t)+\int_{0}^{t}(1+s)^{2\dg}\frac{|\rep'(s)|^{2}}{\cep(s)}\,ds
	\leq k_{24}\ep^{2}(1+t)^{\dg}\lambda^{2}(t)+
	\frac{1}{2}\sup_{\tau\in[0,T]}\Eep(\tau)
	\quad\quad
	\forall t\in[0,T].
	\label{est:pre-sup}
\end{equation}

Let us forget for a while the integral in the left-hand side, and let
us take the supremum of both sides for $t\in[0,T]$. Due to the 
monotonicity of $\lambda(t)$ we obtain that
$$\sup_{\tau\in[0,T]}\Eep(\tau) \leq
2k_{24}\ep^{2}(1+T)^{\dg}\lambda^{2}(T).$$

Coming back to (\ref{est:pre-sup}) we deduce now that
$$\int_{0}^{T}(1+s)^{2\dg}\frac{|\rep'(s)|^{2}}{\cep(s)}\,ds
\leq 2k_{24}\ep^{2}(1+T)^{\dg}\lambda^{2}(T).$$

Since $T$ is arbitrary, we have actually proved that
\begin{equation}
	(1+t)^{2\dg}\left(
	\ep\frac{|\rep'(t)|^{2}}{\cep(t)}+|A^{1/2}\roep(t)|^{2}\right)+
	\int_{0}^{t}(1+s)^{2\dg}\frac{|\rep'(s)|^{2}}{\cep(s)}\,ds \leq
	k_{25}\ep^{2}(1+t)^{\dg}\lambda^{2}(t)
	\label{est:int-2a}
\end{equation}
for every $t\geq 0$.  Plugging (\ref{est:int-2a}) into
(\ref{est:lin-1}) we obtain also that
\begin{equation}
	\int_{0}^{t}(1+s)^{2\dg-1}|A^{1/2}\roep(s)|^{2}\,ds\leq
	k_{26}\ep^{2}(1+t)^{\dg}\lambda^{2}(t)
	\quad\quad
	\forall\, t\geq 0.
	\label{est:int-2b}
\end{equation}

Exploiting once again (\ref{hp:cep-1}), all the estimates in
statement~(1) of Theorem~\ref{thm:main} follow from (\ref{est:int-2a})
and (\ref{est:int-2b}).

\subparagraph{\textmd{\emph{Third energy estimate}}}

We prove that $\roep(t)$ and $\rep(t)$ satisfy the estimates in
statement~(2) of Theorem~\ref{thm:main}.  From the previous step we
already know that
\begin{equation}
	|A^{1/2}\roep(t)|^{2}\leq
	k_{27}\ep^{2}\frac{\lambda^{2}(t)}{(1+t)^{\dg}}
	\quad\quad
	\forall t\geq 0.
	\label{hp:auq-roep}
\end{equation}

Now we observe that, since equations (\ref{eqn:roep-lin}) and
(\ref{eqn:rep-lin}) are linear, we have that the functions
$A^{1/2}\roep(t)$ and $A^{1/2}\rep(t)$ are solutions of analogous
equations, just with $A^{1/2}\gep(t)$ instead of $\gep(t)$.  In this
context (\ref{hp:auq-roep}) and (\ref{hp:gep-2}) play the role of
(\ref{th:roep}) and (\ref{hp:gep-1}), respectively.

Therefore, applying the previous estimates to these new equations, we
immediately obtain that
\begin{equation}
	|A\roep(t)|^{2}\leq
	k_{28}\ep^{2}\frac{\lambda^{2}(t)}{(1+t)^{\dg}}
	\quad\quad
	\forall t\geq 0,
	\label{est:pre-remaining}
\end{equation}
and the integral estimate in statement~(2) of Theorem~\ref{thm:main}.
It remains to prove that
\begin{equation}
	|\rep'(t)|^{2}\leq k_{29}\ep^{2}\frac{\lambda^{2}(t)}{(1+t)^{\dg+2}}
	\quad\quad
	\forall t\geq 0.
	\label{est:remaining}
\end{equation}

To this end, we consider the energy
$$\Gep(t):=(1+t)^{\dg}\frac{|\rep'(t)|^{2}}{\cep^{2}(t)}.$$

Exploiting (\ref{eqn:rep-lin}), with some computations we obtain that
\begin{eqnarray}
	\Gep'(t) & = & -\frac{1}{\ep}(1+t)^{\dg}\left(
	2+2\ep\frac{\cep'(t)}{\cep(t)}-
	\frac{\dg\ep}{1+t}\right)
	\frac{|\rep'(t)|^{2}}{\cep^{2}(t)} 
	\nonumber \\
	\noalign{\vspace{1ex}}
	 &  & -\frac{2}{\ep}\frac{(1+t)^{\dg}}{\cep(t)}
	 \langle A\roep(t),\rep'(t)\rangle
	 +\frac{2}{\ep}\frac{(1+t)^{\dg}}{\cep^{2}(t)}
	 \langle\gep(t),\rep'(t)\rangle.
	 \label{deriv:Gep}
\end{eqnarray}

Let us estimate the three terms. As for the first one, from 
(\ref{hp:cep-2}) we have that
\begin{equation}
	2+2\ep\frac{\cep'(t)}{\cep(t)}- \frac{\dg\ep}{1+t}\geq 1
	\quad\quad
	\forall t\geq 0
	\label{est:Gep-1}
\end{equation}
provided that $\ep$ is small enough. Let us consider now the second 
term. From (\ref{est:pre-remaining}) we have that
\begin{eqnarray}
	-\frac{(1+t)^{\dg}}{\cep(t)}\cdot\langle A\roep(t),\rep'(t)\rangle &
	\leq & (1+t)^{\dg}\cdot\frac{|\rep'(t)|}{\cep(t)}\cdot|A\roep(t)|
	\nonumber  \\
	 & = & \sqrt{\Gep(t)}\cdot(1+t)^{\dg/2}\cdot|A\roep(t)|
	\nonumber  \\
	\noalign{\vspace{1ex}}
	 & \leq & k_{30}\sqrt{\Gep(t)}\cdot\ep\lambda(t).
	\label{est:Gep-2}
\end{eqnarray}

As for the third term, we exploit the estimate from below in
(\ref{hp:cep-1}), and our assumption~(\ref{hp:gep-3}). We obtain that
\begin{eqnarray}
	\frac{(1+t)^{\dg}}{\cep^{2}(t)}\cdot\langle\gep(t),\rep'(t)\rangle & 
	\leq & (1+t)^{\dg}\cdot\frac{|\rep'(t)|}{\cep(t)}\cdot
	\frac{|\gep(t)|}{\cep(t)}
	\nonumber  \\
	 & \leq & \sqrt{\Gep(t)}\cdot(1+t)^{\dg/2}\cdot\frac{1}{\cep(t)}
	 \cdot|\gep(t)|
	\nonumber  \\
	 & \leq & \sqrt{\Gep(t)}\cdot k_{31}(1+t)^{1+\dg/2}\cdot|\gep(t)|
	\nonumber  \\
	\noalign{\vspace{1ex}}
	 & \leq & k_{32}\sqrt{\Gep(t)}\cdot\ep\lambda(t).
	\label{est:Gep-3}
\end{eqnarray}

Plugging (\ref{est:Gep-1}) through (\ref{est:Gep-3}) 
into (\ref{deriv:Gep}) we obtain that
$$\Gep'(t)\leq -\frac{1}{\ep}\Gep(t)+
\frac{1}{\ep}\sqrt{\Gep(t)}\cdot k_{33}\ep\lambda(t)=
-\frac{1}{\ep}\sqrt{\Gep(t)}\left(
\sqrt{\Gep(t)}-k_{33}\ep\lambda(t)\right).$$

Since $\Gep(0)=0$, from Lemma~\ref{lemma:ODE-solito} we conclude that
$\Gep(t)\leq k_{33}^{2}\ep^{2}\lambda^{2}(t)$ for every $t\geq 0$. 
Thanks to the estimate from above in~(\ref{hp:cep-1}), this
is equivalent to (\ref{est:remaining}).\qed

\subsection{Conclusion of proof of Theorem~\ref{thm:main}}

Let $\cep(t)$ and $c(t)$ be defined by (\ref{defn:c-cep}).  Let
$\gep(t)$ be defined by (\ref{defn:gep}).  We already know that
(\ref{th:roep}) holds true.  Thanks to Proposition~\ref{prop:linear},
it is enough to show that assumptions (\ref{hp:cep-1}) through
(\ref{est:int-2}) are satisfied.

Assumption (\ref{est:int-2}) is trivial when $\lambda(t)\equiv 1$, 
and follows from a simple integration by parts when 
$\lambda(t)=\log(e+t)$.

Let us consider the assumptions on $\cep(t)$.
Estimate (\ref{hp:cep-1}) immediately follows from (\ref{est:auq}).
Moreover, since
$$\frac{|\cep'(t)|}{\cep(t)}=
\frac{2\gamma\left|
\langle A^{1/2}\uep(t),A^{1/2}\uep'(t)\rangle\right|}{|A^{1/2}\uep(t)|^{2}}
\leq 2\gamma\frac{|\uep'(t)|\cdot|A\uep(t)|}{|A^{1/2}\uep(t)|^{2}},$$
estimate (\ref{hp:cep-2}) follows from (\ref{est:auq}) through 
(\ref{est:u'}) (in this point we need the estimate from below for 
$|A^{1/2}\uep(t)|$).

In order to prove estimates on $\gep(t)$, we first estimate
$\cep(t)-c(t)$.  To this end, we apply the mean value theorem to the
function $\sigma^{\gamma}$, and we obtain the inequality
$$|y^{\gamma}-x^{\gamma}|\leq\gamma
\max\{y^{\gamma-1},x^{\gamma-1}\}\cdot|y-x| \quad\quad
\forall x\geq 0,\ \forall y\geq 0.$$

Setting $y:=|A^{1/2}\uep(t)|^{2}$ and 
$x:=|A^{1/2}u(t)|^{2}$, it follows that
\begin{equation}
	|\cep(t)-c(t)|\leq\gamma
	\max\left\{|A^{1/2}\uep|^{2(\gamma-1)},
	|A^{1/2}u|^{2(\gamma-1)}\right\}\cdot
	\left||A^{1/2}\uep|^{2}-|A^{1/2}u|^{2}\right|.
	\label{est:cep-c-1}
\end{equation}

From (\ref{est:auq}) and (\ref{th:dec-par}) with $j=1$ we have that
\begin{equation}
	\max\left\{|A^{1/2}\uep(t)|^{2(\gamma-1)},
	|A^{1/2}u(t)|^{2(\gamma-1)}\right\}\leq
	\frac{k_{1}}{(1+t)^{1-1/\gamma}}.
	\label{est:cep-c-2}
\end{equation}

Moreover, arguing as in (\ref{est:R-3}), we obtain that
\begin{equation}
	\left||A^{1/2}\uep(t)|^{2}-|A^{1/2}u(t)|^{2}\right|\leq
	k_{2}\frac{|\roep(t)|}{(1+t)^{1/(2\gamma)}}.
	\label{est:cep-c-3}
\end{equation}

From (\ref{est:cep-c-1}), (\ref{est:cep-c-2}), (\ref{est:cep-c-3}), 
and (\ref{th:roep}) we conclude that
$$|\cep(t)-c(t)|\leq k_{3}
\frac{\ep\lambda(t)}{(1+t)^{1+(\delta-1)/(2\gamma)}}.$$

From (\ref{th:dec-par}) with $j=2$ we have therefore that
\begin{eqnarray}
	|\gep(t)|^{2} & \leq & 2(\cep(t)-c(t))^{2}|Au(t)|^{2}+
	2\ep^{2}|u''(t)|^{2}
	\nonumber  \\
	 & \leq & 
	 k_{4}\frac{\ep^{2}\lambda^{2}(t)}{(1+t)^{2+\delta/\gamma}}+
	 2\ep^{2}|u''(t)|^{2}.
	\label{est:gep}
\end{eqnarray}

At this point (\ref{hp:gep-1}) follows from (\ref{th:par-v''}) and
(\ref{est:int-2}).
Moreover (\ref{hp:gep-2}) follows in an analogous way exploiting
(\ref{th:par-auq-v''}) instead of (\ref{th:par-v''}).

Finally, from (\ref{est:gep}) and (\ref{th:par-v''-point}) we obtain 
that
$$|\gep(t)|^{2}\leq k_{4}
\frac{\ep^{2}\lambda^{2}(t)}{(1+t)^{2+\delta/\gamma}}+
k_{5}\frac{\ep^{2}}{(1+t)^{4+1/\gamma}}\leq
k_{6}\frac{\ep^{2}\lambda^{2}(t)}{(1+t)^{2+\delta/\gamma}},$$
where in the last inequality we used that $2+\delta/\gamma\leq 
4+1/\gamma$. This proves~(\ref{hp:gep-3}), and completes the proof of 
Theorem~\ref{thm:main}.\qed

\setcounter{equation}{0}
\section{Open problems}\label{sec:open}

The main open problems in the theory of dissipative Kirchhoff
equations have been stated in the last section
on~\cite{gg:survey-diss}. In particular, this paper gives a partial 
answer to the sixth problem presented therein.

Here we state some open questions which are more closely related to 
the specific degenerate nonlinearity considered in this paper.
\begin{itemize}
	\item \emph{Open problem 1}.  In the case where $\delta_{0}\geq
	2\gamma+1$, is the term $\lambda(t)$ really needed in the
	estimates of Theorem~\ref{thm:main}?

	\item \emph{Open problem 2}.  Determine the better decay-error
	estimates which are true without the $(\nu,\delta_{0})$-assumption
	on initial data.  We suspect that nothing more than
	(\ref{est:de-trivial}) can be true for general data.

	\item \emph{Open problem 3}.  Is it possible to extend the theory
	to the case $\gamma\in(0,1)$?  

	\item  \emph{Open problem 4}.  Is it possible to extend the 
	theory to weak dissipation terms of the form $(1+t)^{-p}\uep'(t)$ 
	with $p\leq 1$?
\end{itemize}

\label{NumeroPagine}

\end{document}